\documentclass[12pt,letterpaper]{article}%
\usepackage{graphics}
\usepackage{color}
\usepackage{amsmath}
\usepackage{amsthm}
\usepackage{amssymb}
\usepackage[english]{babel}
\usepackage{amsfonts}
\usepackage{graphicx}%
\newtheorem{theorem}{Theorem}[section]
\numberwithin{equation}{section}

\newtheorem{corollary}[theorem]{Corollary}

\newtheorem{hypothesis}[theorem]{Hypothesis}
\newtheorem{lemma}[theorem]{Lemma}

\newtheorem{proposition}[theorem]{Proposition}

\newtheorem{remark}[theorem]{Remark}

\DeclareMathOperator*{\coloneq}{\ensuremath{\mbox{: =}}}

\begin{document}

\title{Superdiffusivity for a Brownian polymer in a continuous Gaussian environment}
\author{S\'ergio Bezerra \thanks{This author's research partially supported by CAPES.}
\hspace{0.3cm} Samy Tindel\\{\small \textit{Institut Elie Cartan, Universit\'e de Nancy 1}}\\{\small \textit{BP 239, 54506-Vandoeuvre-l\`es-Nancy, France}}\\{\small \texttt{[bezerra,tindel]@iecn.u-nancy.fr}}\vspace*{0.05in}
\and Frederi Viens \thanks{This author's research partially supported by NSF grant
no. : 0204999.}\\{\small \textit{Dept. Statistics \& Dept. Mathematics, Purdue University}}\\{\small \textit{150 N. University St., West Lafayette, IN 47907-2067, USA}}\\{\small \texttt{viens@purdue.edu}}\vspace*{0.05in}}
\date{March 1, 2007}
\maketitle

\begin{abstract}
This paper provides information about the asymptotic behavior of a
one-dimen\-sio\-nal Brownian polymer in random medium represented by a
Gaussian field $W$ on ${\mathbb{R}}_{+}\times{\mathbb{R}}$ assumed to be white
noise in time and function-valued in space. According to the behavior of the
spatial covariance of $W$, we give a lower bound on the power growth
(wandering exponent) of the polymer when the time parameter goes to infinity:
the polymer is proved to be superdiffusive, with a wandering exponent
exceeding any $\alpha<3/5$.

\end{abstract}

\vspace{1cm}

\vspace{1cm}


\vspace{1cm}


\newpage

\section{Introduction}

This paper is concerned with a model for a one-dimensional directed Brownian
polymer in a Gaussian random environment (random medium) which can be briefly
described as follows: the polymer itself, in the absence of any random
environment, will simply be modeled by a Brownian motion $b=\{b_{t};t\geq0\}$,
defined on a complete filtered probability space $(\mathcal{C},\mathcal{F}%
,(\mathcal{F}_{t})_{t\geq0},(P_{b}^{x})_{x\in{\mathbb{R}}})$, where $P_{b}%
^{x}$ stands for the Wiener measure starting from the initial condition $x$.
The corresponding expected value will be denoted\footnote{This notation, which
employs a subscript $b$ in a somewhat abusive way to indicate that an average
with respect to the distribution of $b$ is taken, is now common in random
medium theory, and has the advantage of reminding the reader that the
randomness being averaged out is that of the Brownian $b$, not of the medium.}
by $E_{b}^{x}$, or simply by $E_{b}$ when $x=0$. One may assume that
$\mathcal{C}$ is the space of continuous functions started at $0$%
.\vspace{0.3cm}

The random environment will be represented by a centered Gaussian field $W$
indexed by ${\mathbb{R}}_{+}\times{\mathbb{R}}$, defined on another complete
probability space $(\Omega,\mathcal{G},\mbox{{\bf P}})$ independent of $b$'s
space canonical space. Denoting by $\mbox{{\bf E}}$ the expected value with
respect to $\mbox{{\bf P}}$, the covariance structure of $W$ is given by
\begin{equation}
\mbox{{\bf E}}\left[  W(t,x)W(s,y)\right]  =[t\wedge s]\,Q(x-y), \label{E:en2}%
\end{equation}
for a given homogeneous covariance function $Q:{\mathbb{R}}\rightarrow
{\mathbb{R}}$ satisfying some growth conditions that will be specified later
on. In particular, the function $t\mapsto\lbrack Q(0)]^{-1/2}W(t,x)$ is a
standard Brownian motion for any fixed $x\in{\mathbb{R}}$, and for every fixed
$t\in{\mathbb{R}}_{+}$, the process $x\mapsto t^{-1/2}W(t,x)$ is a homogeneous
Gaussian field on ${\mathbb{R}}$ with covariance function $Q$.\vspace{0.3cm}

Once $b$ and $W$ are defined, the polymer measure itself can be described as
follows: for any $t>0$, the energy of a given path (or configuration) $b$ on
$[0,t]$, under the influence of the random environment $W$, is given by the
\emph{Hamiltonian}
\begin{equation}
-H_{t}(b)=\int_{0}^{t}W(ds,b_{s}). \label{nrj}%
\end{equation}
A completely rigorous meaning for this integral will be given in the next
section, but for the moment, notice that for any fixed path $b$, $H_{t}(b)$ is
a centered Gaussian random variable with variance $tQ(0)$. Based on this
Hamiltonian, for any $x\in{\mathbb{R}}$, and a given constant $\beta$
(interpreted as the inverse of the temperature of the system), we define our
(random) polymer measure $G_{t}^{x}$ (with $G_{t}:=G_{t}^{0}$) as follows:
\begin{equation}
dG_{t}^{x}(b)=\frac{e^{-\beta H_{t}(b)}}{Z_{t}^{x}}dP_{b}^{x}(b),\quad
\mbox{ with }\quad Z_{t}^{x}=E_{b}^{x}\left[  e^{-\beta H_{t}(b)}\right]  .
\label{defgt}%
\end{equation}
\vspace{0.3cm}

After early results in the Mathematical Physics literature (see \cite{DSp},
\cite{IS}), links between martingale theory and directed polymers in random
environments were established in \cite{Bo}, \cite{AZ}, and over the last few
years, several papers have shed some light on different types of polymer
models: the case of random walks in discrete potential is treated for instance
in \cite{CH}, the case of Gaussian random walks is in \cite{Me}, \cite{Peter},
and the case of Brownian polymers in a Poisson potential is considered in
\cite{CY2}. On the other hand, the second author of this paper has undertaken
in \cite{RT} the study of the polymer measure $G_{t}$ defined by
(\ref{defgt}). This latter model, which is believed to behave similarly to the
other directed polymers mentioned above, has at least one advantage, from our
point of view: it can be tackled with a wide variety of methods, some of which
are new to the field: scaling invariances for both $b$ and $W$, stochastic
analysis, Gaussian tools. Our long-term goal is to exploit such tools in order
to get a rather complete description of the asymptotic behavior of the measure
$G_{t}$.\vspace{0.3cm}

In the present article, we undertake this task by investigating the so-called
\emph{wandering exponent} $\alpha$
, which measures the growth of the polymer when $t$ tends to $\infty$, and can
be defined informally by the fact that, under the measure $G_{t}$,
$\sup_{s\leq t}|b_{s}|$ behave like $t^{\alpha}$ for large times $t$. This
kind of exponent has been studied in different contexts in \cite{CY2},
\cite{Me}, \cite{Peter}, \cite{Pi} and \cite{Wu}, yielding the conclusion
that, for a wide number of models in dimension one, we should have
$3/5\leq\alpha\leq3/4$. The true exponent conjectured by physicists is
$\alpha=2/3$.

Our understanding, from references \cite{FH}, \cite{HH}, and \cite{KS}, is
that physicists have come to this conjecture in dimension one, based on
simulations (e.g. \cite{HH}) and on theoretical evidence as well as physical
heuristics (in \cite{FH} where $\alpha$ is denoted by $\zeta$). The lower
bound $\alpha\geq3/5$ is confirmed mathematically in partially discrete
settings (e.g. \cite{Peter}). Our section \ref{STRATEGY} provides an
explanation of how our quantitative results confirm that $\alpha$ should be no
less than $3/5$ if the environment's spatial memory, i.e. its spatial
correlation range, is short enough (cubic decay rate), and that
superdiffusivity ($\alpha>1/2$) is only guaranteed if this memory is not too
long (decay rate exponent exceeding $5/2$). These long-spatial-memory
situations are ones which do not seem to be considered in the mathematical or
physical literature, so it is possible that the conjecture $\alpha=2/3$ may
not apply, although at this stage we have no evidence of any example of an
upper bound result implying $\alpha<2/3$.

In this paper, we will see that, for our model, we have $\alpha\geq3/5$. More
specifically we will prove the following.

\begin{theorem}
\label{T:upper} Let $\beta$ be any strictly positive real number. Assume that
$Q:{\mathbb{R}}\rightarrow{\mathbb{R}}$ defined by (\ref{E:en2}) is a
symmetric positive function, decreasing on ${\mathbb{R}}_{+}$ and such that
for some constant $\theta>0$,%
\begin{equation}
Q(x)=O\Bigl(\frac{1}{|x|^{3+\theta}}\Bigr),\quad\mbox{ as }\quad
x\rightarrow\pm\infty, \label{decq}%
\end{equation}
In particular $Q\left(  0\right)  <\infty$, which implies that $W$ defined in
(\ref{E:en2}) is function-valued in $x$. Then, for any $\varepsilon>0$, we
have
\begin{equation}
\lim_{t\rightarrow\infty}\mbox{{\bf P}}\left[  \frac{1}{t^{\frac{3}%
{5}-\varepsilon}}\langle\sup_{s\leq t}|b_{s}|\rangle_{t}\geq1\right]  =1
\label{pgrow}%
\end{equation}
where $\langle\cdot\rangle_{t}$ denotes expectation with respect to the
polymer measure $dG_{t}^{x}\left(  b\right)  $ in (\ref{defgt}).
\end{theorem}

Our proof of this result inspires itself with some of the steps of Peterman's
work in \cite{Peter}, where the same kind of growth bound has been established
for a random walk in a Gaussian potential. Notice that, beyond generalizing
his work from discrete to continuous space, we have been able to extend
Petermann's result to a wider class of environments: indeed we prove the
relation (\ref{pgrow}) holds as soon as $Q$ satisfies the mild correlation
decay assumption (\ref{decq}); Peterman assumed an exponential decay for $Q$.
Moreover, many arguments had to be changed in order to pass from the random
walk to the Brownian case. Having said all this, we must express our debt to
Peterman's work which, unfortunately, has not been published beyond this Ph.D.
dissertation \cite{Peter} as directed by Erwin Bolthausen.

From the physical standpoint, it is worth noting that the above
superdiffusivity theorem (wandering exponent $\alpha>1/2$), which obviously
does not hold for $\beta=0$ (absence of random environment), holds nonetheless
for all $\beta>0$, i.e. all temperatures. This is in contrast to the notion of
strong disorder, defined and described at the end of the next section, a
concept that we will study in detail in a separate publication. However, taken
in a naive and intuitive sense, strong disorder is morally implied by
superdiffusivity; lower bounds on wandering exponents that exceed $1/2$ thus
appear as a convenient quantitative way of measuring this disorder, which is
proved here to hold uniformly for all temperatures.

This paper is structured as follows. Section \ref{PRELIM} defines the random
environment $W$ and the Hamiltonian $H_{t}\left(  b\right)  $ rigorously, and
discusses the relation between our wandering exponent $\alpha$ and the concept
of strong disorder. Section \ref{STRATEGY} discusses the meaning of our main
technical hypothesis \ref{decq}, what happens when one tries to weaken it, and
a related open problem on the interplay between superdiffusivity and random
environment correlation range. Section \ref{STRATEGY} also presents the main
strategy for proving Theorem \ref{T:upper}. The remainder of the paper is
devoted to proving this theorem. Section \ref{INITCOV} calculates the
asymptotic correlation structure of space-time averages of $W$. Section
\ref{INTERACT} calculates similar asymptotics describing the interaction
between $b$ and $W$. Section \ref{girsanov} presents an application of
Girsanov's theorem for $b$ which estimates the penalization needed to force
distant portions of $b$ back near the origin. Finally, with all these
quantitative tools in hand, the proof of the theorem is completed in Section
\ref{Lenvproof}, which also contains a detailed heuristic description of this
part of the proof.

The authors of this paper express their thanks to two referees whose detailed
comments resulted in corrections and other improvements over an earlier
version of this paper.

\section{Preliminaries; the partition function; strong disorder\label{PRELIM}}

In this section, we will first recall some basic facts about the partition
function $Z_{t}$, and then give briefly some notions of Gaussian analysis
which will be used later on. Let us recall that $W$ is a centered Gaussian
field defined on ${\mathbb{R}}_{+}\times{\mathbb{R}}$, which can also be seen
as a Gaussian family $\{W(\varphi)\}$ indexed by tests functions
$\varphi:{\mathbb{R}}_{+}\times{\mathbb{R}}\rightarrow{\mathbb{R}}$, where
$W(\varphi)$ stands for the Wiener integral of $\varphi$ with respect to $W$:
\[
W(\varphi)=\int_{{\mathbb{R}}}\int_{{\mathbb{R}}_{+}}\varphi(s,x)W(ds,x)dx,
\]
whose covariance structure is given by
\begin{equation}
\mbox{{\bf E}}\left[  W(\varphi)W(\psi)\right]  =\int_{{\mathbb{R}}_{+}%
}\left(  \int_{{\mathbb{R}}\times{\mathbb{R}}}\varphi(s,x)Q(x-y)\psi
(s,y)dxdy\right)  ds, \label{covw}%
\end{equation}
for two arbitrary test functions $\varphi,\psi$.

Let us start here by defining more rigorously the quantity $H_{t}(b)$ given by
(\ref{nrj}), which can be done through a Fourier transform procedure: there
exists (see e.g. \cite{CV} for further details) a centered Gaussian
independently scattered $\mathbb{C}$-valued measure $\nu$ on ${\mathbb{R}}%
_{+}\times{\mathbb{R}}$ such that
\begin{equation}
W(t,x)=\int_{{\mathbb{R}}_{+}\times{\mathbb{R}}}\mathbf{1}_{[0,t]}%
(s)e^{iux}\nu(ds,du). \label{E:rep_int}%
\end{equation}
For every test function $f:{\mathbb{R}}_{+}\times{\mathbb{R}}\rightarrow
\mathbb{C}$, set now
\begin{equation}
\nu(f)\equiv\int_{{\mathbb{R}}_{+}\times{\mathbb{R}}}f(s,u)\nu(ds,du).
\label{defnuf}%
\end{equation}
While the random variable $\nu\left(  f\right)  $ may be complex-valued, to
ensure that it is real valued, it is sufficient to assume that $f$ is of the
form $f\left(  s,u\right)  =f_{1}\left(  s\right)  e^{iuf_{2}\left(  s\right)
}$ for real valued functions $f_{1}$ and $f_{2}$. Then the law of $\nu$ is
defined by the following covariance structure: for any such test functions
$f,g:{\mathbb{R}}_{+}\times{\mathbb{R}}\rightarrow\mathbb{C}$, we have
\begin{equation}
\mbox{{\bf E}}\left[  \nu(f)\nu(g)\right]  =\int_{{\mathbb{R}}_{+}%
\times{\mathbb{R}}}f(s,u)\overline{g(s,u)}\hat{Q}(du)ds, \label{E:cov_X}%
\end{equation}
where the finite positive measure $\hat{Q}$ is the Fourier transform of $Q$
(see \cite{tv99} for details).

From (\ref{E:rep_int}), we see that the It\^{o}-stochastic differential of $W$
in time can be understood as $W\left(  ds,x\right)  :=\int_{u\in{\mathbb{R}}%
}e^{iux}\nu(ds,du)$, or even, if the measure $\hat{Q}\left(  du\right)  $ has
a density $f\left(  u\right)  $ with respect to the Lebesgue measure, which is
typical, as
\[
W\left(  ds,x\right)  :=\int_{u\in{\mathbb{R}}}e^{iux}\sqrt{f\left(  u\right)
}M(ds,du)
\]
where $M$ is a white-noise measure on ${\mathbb{R}}_{+}\times{\mathbb{R}}$,
i.e. a centered independently scattered Gaussian measure with covariance given
by $\mathbf{E}\left[  M\left(  A\right)  M\left(  B\right)  \right]
=m_{Leb}\left(  A\cap B\right)  $ where $m_{Leb}$ is Lebegue's measure on
${\mathbb{R}}_{+}\times{\mathbb{R}}$.

We can go back now to the definition of $H_{t}(b)$: invoking the
representation (\ref{E:rep_int}), we can write
\begin{equation}
-H_{t}(b)=\int_{0}^{t}W(ds,b_{s})=\int_{0}^{t}\int_{{\mathbb{R}}}e^{iub_{s}%
}\nu(ds,du), \label{E:rep_intH}%
\end{equation}
and it can be shown (see \cite{CV}) that the right hand side of the above
relation is well defined for any H\"{o}lder continuous path $b$, by a $L^{2}%
$-limit procedure. Such a limiting procedure can be adapted to the specific
case of constructing $H_{t}\left(  b\right)  $, using the natural time
evolution structure; we will not comment on this further. However, the reader
will surmise that the following remark, give for the sake of illustration, can
be useful: when $\hat{Q}$ has a density $f$, we obtain%
\[
-H_{t}(b)=\iint_{[0,t]\times{\mathbb{R}}}e^{iub_{s}}\sqrt{f\left(  u\right)
}M\left(  ds,du\right)  .
\]

\vspace{0.3cm}

With the so-called \emph{partition function }$Z_{t}^{x}$ defined earlier as
$Z_{t}^{x}=E_{b}\left[  e^{-\beta H_{t}(b)}\right]  $, set
\begin{equation}
p_{t}(\beta):=\frac{1}{t}\mbox{{\bf E}}\left[  \log\left(  Z_{t}^{x}\right)
\right]  , \label{fren}%
\end{equation}
usually called the free energy of the system. By spatial homogeneity of $W$,
$p_{t}(\beta)$ is independent of the initial condition $x\in{\mathbb{R}}$, and
the same holds for the law of $b-x$ under $G_{t}^{x}$, thus without loss of
generality we set $x=0$, hence the notation $E_{b},Z_{t}$,... standing for
$E_{b}^{0},Z_{t}^{0}$, etc. It was shown in \cite{RT} that $\lim
_{t\rightarrow\infty}p_{t}(\beta)=\sup_{t\geq0}p_{t}(\beta)$ exists and is
positive, and that $\mbox{{\bf P}}$-almost surely, $\frac{1}{t}\log Z_{t}$
converges to the same limit$.$ The trivial bound%
\begin{equation}
p(\beta):=\lim_{t\rightarrow\infty}p_{t}(\beta)\leq\frac{\beta^{2}}{2}Q(0)
\label{roughbnd}%
\end{equation}
always holds, but the polymer is said to be in the \emph{strong disorder}
regime if $\lim_{t\rightarrow\infty}\frac{1}{t}\log Z_{t}<\frac{\beta^{2}}%
{2}Q(0)$, which is therefore equivalent to saying that inequality
(\ref{roughbnd}) above is strict. We will show in a separate publication that,
for all $\beta\geq\beta_{0}$, while $p\left(  \beta\right)  \geq c\beta^{4/3}$
for all non-trivial random media $W$ and some constant depending on $W$'s law,
we have the specific strong disorder upper bound $p\left(  \beta\right)  \leq
c\beta^{2-2H/(2H+1)}$ where $H$ is a spatial H\"{o}lder exponent for $W$. Yet
we do not know if these results can be made to hold for small $\beta$. One
would prefer not having any condition on the temperature scale, and physicists
expect strong disorder in our one-dimensional setting for all $\beta>0$, which
is only confirmed mathematically in some cases, such as in \cite{CH} and
\cite{CY2}.

This is where the polymer's superdiffusivity (wandering exponent $\alpha>1/2$)
can be useful to our fully continuous situation. Since the concept of
\textquotedblleft strong disorder\textquotedblright\ was introduced in order
to determine whether the random environment has any significant influence on
polymer paths $b$, it is generally acceptable to say that a polymer with
super-diffusive behavior exhibits \textquotedblleft strong
disorder\textquotedblright. Even though this second definition does not match
the common one given above ($p\left(  \beta\right)  =Q\left(  0\right)
\beta^{2}/2$), it is useful to note that the results of the next section imply
the following (see Corollary \ref{disorder}): if $W$ exhibits decorrelation
that is not too slow, specifically if for large $x$, $Q\left(  x\right)  \leq
cx^{-5/2-\vartheta}$ where $\vartheta>0$, then the polymer is superdiffusive
with exponent any $\alpha<\min\left\{  \frac{1}{2}+\frac{\vartheta
}{6-2\vartheta};3/5\right\}  $, and this form of strong disorder holds for all
$\beta>0$. The specific order of decorrelation $x^{-5/2-\vartheta}\ll
x^{-5/2}$ can be quantified by saying that $W$'s decorrelation is certainly
faster than the well-known order $x^{-2+2H}$ for the increments of fractional
Brownian motion, but the class of such $W$'s still qualifies as containing
long-range correlations (polynomial with moderate power).

We also plan to investigate, in a separate publication, situations in which we
can show the complementary story: we plan to prove that if \emph{weak
disorder} holds, i.e. if $\lim_{t\rightarrow\infty}\frac{1}{t}\log Z_{t}%
=\frac{\beta^{2}}{2}Q(0)$, then the polymer is diffusive, i.e. $\alpha=1/2$.

\section{Discussion of hypothesis and results; strategy of
proof\label{STRATEGY}}

Recall our goal: we will prove that for the polymer measure $G_{t}=G_{t}^{0}$
in (\ref{defgt}), Theorem \ref{T:upper} holds. This theorem gives an
indication of the asymptotic speed of our polymer. Indeed, if we could write
that $\sup_{s\leq t}|b_{s}|\sim t^{\alpha}$ under $G_{t}$ as $t\rightarrow
\infty$, then Theorem \ref{T:upper} would state that the wandering exponent
$\alpha\ $is no smaller than $3/5$. As stated in the introduction, our basic
technical assumption to prove the theorem is the following.

\begin{hypothesis}
\label{env1} We assume that $Q:{\mathbb{R}}\rightarrow{\mathbb{R}}$ defined by
(\ref{E:en2}) is a symmetric positive function, decreasing on ${\mathbb{R}%
}_{+}$ and such that there exists a strictly positive constant $\theta$ such
that
\[
Q(x)=O\Bigl(\frac{1}{|x|^{3+\theta}}\Bigr),\quad\mbox{ as }\quad
x\rightarrow\pm\infty.
\]

\end{hypothesis}

The rate $3+\theta$ can be quantified physically by saying that $W$
decorrelates in space faster than the well-known order $x^{-2+2H}$ for the
increments of fractional Brownian motion with Hurst parameter $H\in(0,1)$, but
the class of $W$'s defined by Hypothesis \ref{env1} still qualifies as
containing long-range correlated noises (polynomial rate with moderate power),
as opposed to exponential correlation decay, found for instance in finite
memory ARCH/GARCH models, and even more so in opposition to the case of
spatial white noise.

The specific correlation decay rate of $Q$ in the above hypothesis appears to
be important in order to obtain the highest possible superdiffusion wandering
exponent $\alpha$ using our technique (any $\alpha<3/5$). The end of Section
\ref{Lenvproof} shows that if one tries use a smaller decay power than
$3+\theta$ above, the result is impeded: $\alpha$ cannot be chosen arbitrarily
close to $3/5$. In Corollary \ref{disorder} and its preceeding discussion, we
prove that if $Q\left(  x\right)  =O\left(  \left\vert x\right\vert
^{-r}\right)  $ with $r\in(5/2,3]$, then we can only guarantee being able to
take $1/2<\alpha<3/\left(  11-2r\right)  $, so superdiffusivity is still
proved, but $\alpha$ arbitrarily close to $3/5$ is disallowed.

Corollary \ref{disorder} thus opens the interesting question of whether, in
continuous space, the Brownian polymer in a Gaussian environment has a
super-diffusive behavior with a wandering exponent determined by the
environment's range/rate of spatial correlations. We do not believe that any
physical conjecture in which $\alpha=2/3$ specifically argues that this should
hold in our continuous space setting. There are other examples in which
scaling limits depend heavily on whether one is in discrete or continuous
space: for instance, in the regime of small diffusion constant (resp.
viscosity) $\kappa$, the almost-sure Lyapunov exponent for the partition
function $Z_{t}$ (resp. Anderson model) is known to depend heavily on the
spatial regularity of $W$ in continuous space (see \cite{FV}), but is know to
be universally of order $1/\log\left(  \kappa^{-1}\right)  $ in discrete space
(see \cite{CKM}). We will not discuss this point further herein.

\begin{remark}
\label{R:nondeg}Hypothesis \ref{env1} immediately implies that $Q\left(
0\right)  <\infty$. Since $\max\left\vert Q\right\vert =Q\left(  0\right)  $
and $Q$ has an integrable tail, we get $Q\in L^{1}({\mathbb{R}})$.

Without loss of generality, we will assume throughout that $Q$ is normalized
so that $\int_{{\mathbb{R}}}Q(x)dx=1$.

The integrability of $Q$ represents a kind of non-degeneracy condition, which
says that the decorrelation of $W$ at distinct sites is not immediate.\vspace
{0.3cm}
\end{remark}

\begin{proof}
[Strategy of the proof for Theorem \ref{T:upper}]For $t,\epsilon>0$, set
\[
A_{t,\epsilon}=\left\{  \mbox{there exists }s_{0}\in\left[  t/2,t\right]
\text{ such that}|b_{s_{0}}|\geq t^{\frac{3}{5}-\frac{\epsilon}{2}}\right\}
.
\]
Then we can write
\begin{align}
\frac{\langle\sup_{s\leq t}|b_{s}|\rangle_{t}}{t^{\frac{3}{5}-\epsilon}}  &
\geq\frac{t^{\frac{\epsilon}{2}}}{t^{\frac{3}{5}-\frac{\epsilon}{2}}%
}\left\langle \sup_{s\leq t}|b_{s}|\,\mathbf{1}_{A_{t,\epsilon}}\right\rangle
_{t}\nonumber\\
&  \geq t^{\frac{\epsilon}{2}}G_{t}\Bigl(A_{t,\epsilon}\Bigr),\nonumber
\end{align}
since $\sup_{s\leq t}|b_{s}|\geq t^{\frac{3}{5}-\frac{\epsilon}{2}}$ on
$A_{t,\epsilon}$. Thus
\begin{equation}
\frac{\langle\sup_{s\leq t}|b_{s}|\rangle_{t}}{t^{\frac{3}{5}-\epsilon}}\geq
t^{\frac{\epsilon}{2}}\left(  1-G_{t}\left(  A_{t,\epsilon}^{c}\right)
\right)  , \label{E:do}%
\end{equation}
where $A_{t,\epsilon}^{c}=\{b;\,\sup_{s\in\lbrack t/2,t]}|b_{s}|\leq
t^{\frac{3}{5}-\frac{\epsilon}{2}}\}$ is the complement of $A_{t,\epsilon}$.
We will start now a discretization procedure in space: for an arbitrary
integer $k$, and $\alpha>0$, set
\[
I_{k}^{\alpha}=t^{\alpha}[2k-1,2k+1),\quad\mbox{ and }\quad L_{k}^{\alpha
}=\left\{  b;\,b_{s}\in I_{k}^{\alpha}\mbox{ for all }s\in\lbrack
t/2,t]\right\}  .
\]
Then $\tilde{A}_{t,\epsilon}=L_{0}^{3/5-\epsilon/2}$, and equation~(\ref{E:do}%
) can be rewritten as
\[
\frac{\langle\sup_{s\leq t}|b_{s}|\rangle_{t}}{t^{\frac{3}{5}-\epsilon}}\geq
t^{\frac{\epsilon}{2}}\left(  1-G_{t}\left(  L_{0}^{\frac{3}{5}-\frac
{\epsilon}{2}}\right)  \right)  .
\]
Set now
\[
Z_{t}^{\alpha}(k)\coloneq E_{b}\left[  \mathbf{1}_{L_{k}^{\alpha}}\exp\left(
-\beta H_{t}(b)\right)  \right]  .
\]
We have
\[
\frac{\langle\sup_{s\leq t}|b_{s}|\rangle_{t}}{t^{\frac{3}{5}-\epsilon}}\geq
t^{\frac{\epsilon}{2}}\left(  1-\frac{Z_{t}^{\frac{3}{5}-\frac{\epsilon}{2}%
}(0)}{E_{b}\left[  \exp\left(  -\beta H_{t}(b)\right)  \right]  }\right)  ,
\]
by definition of $G_{t}$. On the other hand, since the events $L_{k}^{\alpha}$
are disjoint sets we have
\[
E_{b}\left[  \exp\left(  -\beta H_{t}(b)\right)  \right]  \geq\sum
_{k\in\mathbb{Z}}Z_{t}^{\frac{3}{5}-\frac{\epsilon}{2}}(k).
\]
Therefore, we have established that
\begin{equation}
\frac{\langle\sup_{s\leq t}|b_{s}|\rangle_{t}}{t^{\frac{3}{5}-\epsilon}}\geq
t^{\frac{\epsilon}{2}}\left(  1-\frac{Z_{t}^{\frac{3}{5}-\frac{\epsilon}{2}%
}(0)}{Z_{t}^{\frac{3}{5}-\frac{\epsilon}{2}}(0)+Z_{t}^{\frac{3}{5}%
-\frac{\epsilon}{2}}(k)}\right)  , \label{lbsupb}%
\end{equation}
for any integer $k\neq0$. Suppose now that $W\in\mathcal{A}_{t}$, where
$\mathcal{A}_{t}$ is defined as
\[
\mathcal{A}_{t}\coloneq\left\{  W\;;\,\mbox{ There exists }k^{\ast}\neq0\text{
such that }Z_{t}^{\alpha}(k^{\ast})>Z_{t}^{\alpha}(0)\right\}  .
\]
Then, choosing $k=k^{\ast}$ in (\ref{lbsupb}), it is easily seen that
\[
\frac{\langle\sup_{s\leq t}|b_{s}|\rangle_{t}}{t^{\frac{3}{5}-\epsilon}}\geq
t^{\frac{\epsilon}{2}}\Bigl(1-\frac{1}{2}\Bigr)\geq1,
\]
whenever $t$ is large enough. The proof is now easily finished if we can prove
the following lemma:

\begin{lemma}
\label{L:env} Given a positive real number $\alpha\in(1/2;3/5)$ and an
environment $W$ satisfying Hypothesis \ref{env1}, then
\begin{align}
\liminf_{t\rightarrow\infty}\mbox{{\bf P}}( \mathcal{A}_{t})=1.
\end{align}

\end{lemma}

The remainder of this article will now be devoted to the proof Lemma
\ref{L:env}.
\end{proof}

\section{Initial covariance computations\label{INITCOV}}

In order to prove Lemma~\ref{L:env}, we shall begin with a series of
preliminary results, the first of which is a covariance computation, including
precise asymptotic estimations in large time, for space-time averages of the
random environment $W$.

For a given $k\in{\mathbb{Z}}$ and $\alpha>0$, recall that $I_{k}%
:=I_{k}^{\alpha}=t^{\alpha}[2k-1,2k+1)$, and set
\begin{equation}
\tilde{\eta}_{k}=\tilde{\eta}_{k}^{\alpha}:=\frac{1}{t^{(\alpha+1)/2}}%
\int_{\frac{t}{2}}^{t}\int_{I_{k}}W(ds,x)dx. \label{D:eta}%
\end{equation}
Then $\{\tilde{\eta}_{k};k\in{\mathbb{Z}}\}$ is a centered Gaussian vector,
whose covariance matrix will be called $C(t)=(C_{\ell,k}(t))_{\ell
,k\in\mathbb{Z}}$, where
\begin{equation}
C_{\ell,k}(t)=\mbox{{\bf E}}\left[  \tilde{\eta}_{\ell}\tilde{\eta}%
_{k}\right]  =\mbox{\bf Cov}\left(  \tilde{\eta}_{\ell};\,\tilde{\eta}%
_{k}\right)  =\frac{1}{2t^{\alpha}}\int_{I_{k}}\int_{I_{\ell}}Q\left(
x-y\right)  dxdy. \label{covet}%
\end{equation}
where the last equality above follows directly from the definition of $W$'s
covariance in (\ref{covw}).

Here and below, we omit the superscripts $\alpha$ on quantities like
$\tilde{\eta}_{k}^{\alpha}$, $I_{k}^{\alpha}$, $L_{k}^{\alpha}$, etc... We now
proceed to estimate the matrix $C\left(  t\right)  $, and show in particular
that $\lim_{t\rightarrow\infty}C(t)=\mbox{Id}$. This can be interpreted as
saying that the amount of decorrelation of the potential at distant locations
implied by Hypothesis \ref{env1}, is enough to guarantee independence of the
$\tilde{\eta}_{k}$ asymptotically.

\begin{proposition}
\label{P:ord_t} Let $\theta$ be the strictly positive constant defined in
Hypothesis \ref{env1}, and consider $k\in{\mathbb{Z}}$, $\alpha>0$ and
$\tau<\theta\wedge1$. Set also
\[
\lambda:=\frac{1}{C_{0,0}(t)}=\frac{1}{C_{k,k}(t)},
\]
where $C(t)$ has been defined at (\ref{covet}). Then, the elements of $C(t)$
satisfy the following.

\begin{itemize}
\item[(i)] \label{item_i} $\lambda=1+ O\left(  \frac{1}{t^{\alpha}}\right)  $.

\item[(ii)] \label{item_ii} $\lambda\sum_{\ell\neq k}|\ell-k|^{\tau} |
C_{\ell,k}(t)|= O\left(  \frac{1}{t^{\alpha}}\right)  \cdot$
\end{itemize}
\end{proposition}

\begin{proof}
~\vspace*{0.1in}

\noindent\emph{Step 0: initial calculation.} We will only consider the case
$k=0$, the other ones being easily deduced by homogeneity of $W$. Let us first
evaluate $C_{\ell,0}(t)$ for $\ell\geq0$ (here again, the case $\ell<0$ is
similar, since $Q$ is a symmetric function). Then, a direct application of
(\ref{covet}) gives
\[
C_{\ell,0}(t)=\frac{1}{2t^{\alpha}}\int_{t^{\alpha}(2\ell-1)}^{t^{\alpha
}(2\ell+1)}\int_{-t^{\alpha}}^{t^{\alpha}}Q(x-y)dxdy.
\]
Set now
\[
(I):=\frac{1}{2t^{\alpha}}\left[  \int_{t^{\alpha}(2\ell-1)}^{t^{\alpha}%
(2\ell+1)}\int_{-\infty}^{-t^{\alpha}}Q(x-y)dxdy+\int_{t^{\alpha}(2\ell
-1)}^{t^{\alpha}(2\ell+1)}\int_{t^{\alpha}}^{\infty}Q(x-y)dxdy\right]  .
\]
Since $\int_{R}Q(x-y)dx=1$ for any $y\in{\mathbb{R}}$, it is easily checked
that
\begin{equation}
C_{\ell,0}(t)=1-\left(  I\right)  . \label{E:clo}%
\end{equation}
Then, a series of changes of variable yields
\begin{align*}
(I)  &  =\frac{1}{2t^{\alpha}}\left[  \int_{t^{\alpha}(2\ell-1)}^{t^{\alpha
}(2\ell+1)}\int_{-\infty}^{-t^{\alpha}-y}Q(u)dudy+\int_{t^{\alpha}(2\ell
-1)}^{t^{\alpha}(2\ell+1)}\int_{t^{\alpha}-y}^{\infty}Q(u)dudy\right] \\
&  =\frac{1}{2t^{\alpha}}\left[  \int_{-t^{\alpha}(2\ell+2)}^{-t^{\alpha
}(2\ell)}\int_{-\infty}^{\hat{z}}Q(u)dud\hat{z}+\int_{-t^{\alpha}(2\ell
)}^{-t^{\alpha}(2\ell-2)}\int_{z}^{\infty}Q(u)dudz\right]  ,
\end{align*}
where we have set $\hat{z}=-t^{\alpha}-y$ and $z=t^{\alpha}-y$. Thus, denoting
by $\bar{F}(z)$ the quantity $\int_{z}^{\infty}Q(u)du$, we get
\begin{align}
(I)  &  =\frac{1}{2t^{\alpha}}\left[  \int_{-t^{\alpha}(2\ell+2)}^{-t^{\alpha
}(2\ell)}\left(  1-\bar{F}(\hat{z})\right)  d\hat{z}+\int_{-t^{\alpha}(2\ell
)}^{-t^{\alpha}(2\ell-2)}\bar{F}\left(  z\right)  )dz\right] \nonumber\\
&  =1-\frac{1}{2t^{\alpha}}\int_{-t^{-\alpha}(2\ell+2)}^{-t^{\alpha}(2\ell
)}\bar{F}(z)dz+\frac{1}{2t^{\alpha}}\int_{-t^{\alpha}(2\ell)}^{-t^{\alpha
}(2\ell-2)}\bar{F}\left(  z\right)  )dz. \label{E:est}%
\end{align}
Putting together~(\ref{E:clo}) and~(\ref{E:est}) one obtains, for any
$\ell\geq0$,
\begin{equation}
C_{\ell,0}(t)=\frac{1}{2t^{\alpha}}\left[  \int_{-t^{\alpha}(2\ell
+2)}^{-t^{\alpha}(2\ell)}\bar{F}(z)dz-\int_{-t^{\alpha}(2\ell)}^{-t^{\alpha
}(2\ell-2)}\bar{F}(z)dz\right]  . \label{E:cl0}%
\end{equation}
\vspace*{0.1in}

\noindent\emph{Step 1: proving item (i).} We are now ready to prove item (i).
By symmetry of $Q$, we have $1-\bar{F}\left(  -z\right)  =\bar{F}\left(
z\right)  $. Thus for $\ell=0$, equation (\ref{E:cl0}) becomes
\begin{equation}
C_{0,0}(t)=\frac{1}{2t^{\alpha}}\left[  \int_{-2t^{\alpha}}^{0}\left(
1-\bar{F}(-z)\right)  dz-\int_{0}^{2t^{\alpha}}\bar{F}(z)dz\right]
=1-\frac{1}{t^{\alpha}}\int_{0}^{2t^{\alpha}}\bar{F}(z)dz.
\end{equation}
Now, using the fact that
\begin{equation}
\bar{F}(z)\leq c\left(  1\wedge|z|^{-(2+\theta)}\right)  , \label{tailbf}%
\end{equation}
which follows directly from Hypothesis \ref{env1}, it is easily seen that
$C_{0,0}(t)=1+O(t^{-\alpha})$, which ends the proof of item (i).\vspace
*{0.1in}

\noindent\emph{Step 2: proving item (ii).} In order to show item (ii), we deal
with $\ell=1$ separately from the other cases. Beginning with $\ell\geq2$, we
first get the obvious derivative $\bar{F}^{\prime}\left(  z\right)  =-Q\left(
z\right)  $, and we will use the fact that $Q$ is decreasing on ${\mathbb{R}%
}_{+}$ to bound this latter function on an interval in ${\mathbb{R}}_{+}$ by
its value at the left endpoint. Invoking the fact that $\bar{F}(-v)=1-\bar
{F}(x)$, we may thus write from equation~(\ref{E:cl0})%
\begin{align*}
|C_{\ell,0}(t)|  &  =\frac{1}{2t^{\alpha}}\left\vert \int_{-t^{\alpha}(2\ell
)}^{-t^{\alpha}(2\ell-2)}\left[  \bar{F}(z-2t^{\alpha})-\bar{F}(z)\right]
dz\right\vert \\
&  =\frac{1}{2t^{\alpha}}\left\vert \int_{-t^{\alpha}(2\ell)}^{-t^{\alpha
}(2\ell-2)}\left[  \bar{F}(-z+2t^{\alpha})-\bar{F}(-z)\right]  dz\right\vert
\\
&  =\frac{1}{2t^{\alpha}}\left\vert \int_{t^{\alpha}(2\ell-2)}^{2t^{\alpha
}\ell}\left[  \bar{F}(z+2t^{\alpha})-\bar{F}(z)\right]  dz\right\vert \\
&  =\frac{1}{2t^{\alpha}}\left\vert \int_{t^{\alpha}(2\ell-2)}^{2t^{\alpha
}\ell}\left(  -\int_{z}^{z+2t^{\alpha}}Q\left(  x\right)  dx\right)
dz\right\vert \\
&  \leq2t^{\alpha}Q\left(  t^{\alpha}(2\ell-2)\right) \\
&  \leq ct^{-\alpha\left(  2+\theta\right)  }\left(  2\ell-2\right)
^{-3-\theta},
\end{align*}
where the last step holds by Hypothesis \ref{env1} for some constant $c>0$. We
immediately obtain%
\begin{align*}
\sum_{\ell=2}^{\infty}|C_{\ell,0}(t)|\ell^{\tau}  &  \leq ct^{-\alpha\left(
2+\theta\right)  }\sum_{\ell=2}^{\infty}\left(  2\ell-2\right)  ^{-3-\theta
}\ell^{\tau}\\
&  \leq cK_{\tau,\theta}t^{-\alpha\left(  2+\theta\right)  }%
\end{align*}
for some constant $K_{\tau,\theta}$ as soon as $\tau<2+\theta$, which is
clearly satisfied by the assumption on $\tau$, and leads to an upper bound in
the series in item (ii) which is amply sufficient to prove the proposition,
except for the term $\ell=1$, with which we deal now.

To finish the proof of the proposition, it is indeed sufficient to prove that
$t^{\alpha}C_{1,0}$ is bounded. We first evaluate this quantity from
(\ref{E:cl0}):%
\begin{align*}
t^{\alpha}C_{1,0}  &  =\int_{-4t^{\alpha}}^{-2t^{\alpha}}\bar{F}%
(z)dz-\int_{-2t^{\alpha}}^{0}\bar{F}(z)dz=\int_{-2t^{\alpha}}^{0}\left(
\bar{F}(z-2t^{\alpha})-\bar{F}(z)\right)  dz\\
&  =\int_{-2t^{\alpha}}^{0}\left(  \int_{z-2t^{\alpha}}^{z}Q\left(  x\right)
dx\right)  dz=\int_{0}^{2t^{\alpha}}\left(  \int_{-z-2t^{\alpha}}^{-z}Q\left(
x\right)  dx\right)  dz\\
&  =\int_{0}^{2t^{\alpha}}\left(  \int_{z}^{z+2t^{\alpha}}Q\left(  x\right)
dx\right)  dz.
\end{align*}
Next we separate the first unit of the $z$-integral from its remainder:
$t^{\alpha}C_{1,0}=A+B$ where we define $A:=\int_{0}^{1}\left(  \int
_{z}^{z+2t^{\alpha}}Q\left(  x\right)  dx\right)  dz$ and $B:=\int
_{1}^{1\wedge2t^{\alpha}}\left(  \int_{z}^{z+2t^{\alpha}}Q\left(  x\right)
dx\right)  dz$. Since $\int_{{\mathbb{R}}}Q=1$, we immediately have $A\leq1$
which is the only term to deal with when $t\leq2^{-1/\alpha}$. When
$t>2^{-1/\alpha}$, for the term $B$, we use Hypothesis \ref{env1}: for some
constant $c$,%
\begin{align*}
B  &  \leq c\int_{1}^{2t^{\alpha}}\left(  \int_{z}^{z+2t^{\alpha}}%
x^{-3-\theta}dx\right)  dz\\
&  =\frac{c}{\left(  \theta+1\right)  \left(  \theta+2\right)  }\left(
1-2^{-\theta}+4^{-\theta-1}\right)  \left(  t^{\alpha}\right)  ^{-\left(
\theta+1\right)  }\leq\frac{c}{\left(  \theta+1\right)  \left(  \theta
+2\right)  }.
\end{align*}
This finishes the proof of the proposition.
\end{proof}

\section{Interaction between $b$ and $W$\label{INTERACT}}

The next step in developping the tools to prove Lemma \ref{L:env} is to get
some quantitative information about the way $b$ interacts with the random
environment $W$ when the Brownian motion is localized by the event $L_{k}$. As
we did with the notation $I_{k}:=I_{k}^{\alpha}$, we are omitting superscripts
$\alpha$ writing only $L_{k}$ instead of $L_{k}^{\alpha}$ from now on.

We begin by introducing two quantities. First, in order to simplify some
$t\,$-dependent normalizers, we renormalize $\tilde{\eta}$ as%
\begin{equation}
\eta_{\ell}:=\frac{t^{\frac{1-\alpha}{2}}}{2}\tilde{\eta}_{\ell}=\frac
{1}{2t^{\alpha}}\int_{\frac{t}{2}}^{t}\int_{I_{\ell}}W(ds,x)dx; \label{D:eta1}%
\end{equation}
we will not need to revert to using $\tilde{\eta}$ in this article. We also
need a vector $v=v(b_{s};t/2\leq s\leq t)$ of ${\mathbb{R}}^{{\mathbb{Z}}}$,
defined for each $\ell\in{\mathbb{Z}}$ by
\begin{equation}
v_{\ell}:=4t^{\alpha-1}\mbox{{\bf E}}\left[  \eta_{\ell}\int_{\frac{t}{2}}%
^{t}W(ds,b_{s})\right]  . \label{defvl}%
\end{equation}
We will prove, in a sense, that $v_{\ell}$ looks like $\mathbf{1}_{\left\{
k\right\}  }\left(  \ell\right)  $ on $L_{k}$. To this end, for a fixed
$k\in{\mathbb{Z}}$, and $\tau<\theta$ (remember that $\theta$ is defined in
Hypothesis \ref{env1}), let us consider the norm $\Vert\cdot\Vert_{\tau,k}$
defined on ${\mathbb{R}}^{{\mathbb{Z}}}$ by
\begin{equation}
\Vert x\Vert_{\tau,k}=|x_{k}|+\sum_{i\neq k}|x_{i}||i-k|^{\tau}.
\label{defntk}%
\end{equation}

\begin{remark}
It will be essential in the sequel to control the decay of $v_{\ell}$, and
also of a quantity $\delta_{\ell}$ (defined later in Proposition
\ref{P:nor_del} as the $\ell$th component of the solution $x$ to the linear
system $C\left(  t\right)  x=v$) when $|\ell|\rightarrow\infty$. It will be
used for instance in relations (\ref{bnputm}) and (\ref{whenQcondchanges}).
This is why we have introduced the norm $\Vert\cdot\Vert_{\tau,k}$ here.
\end{remark}

\subsection{Asymptotics and boundedness of $v$}

We are now ready to state a first result about the interaction between $b$ and
$W$: the behavior of the vector $v$ in large time.

\begin{proposition}
\label{P:nor_v} Suppose $b\in L_{k}$. Then the vector $v$ given by
(\ref{defvl}) satisfies the following properties:

\begin{itemize}
\item[(i)] Let $\|\cdot\|_{\tau,k}$ be the norm defined at (\ref{defntk}).
Then
\[
\| v\|_{\tau,k}- v_{k}= O\left(  \frac{1}{t^{\alpha}}\right)  .
\]

\item[(ii)] For $t$ large enough, there exist two strictly positive real
numbers $\underline{c}$ and $\overline{c}$ such that
\begin{align}
\underline{c}\leq v_{k} \leq\overline{c}.\nonumber
\end{align}

\end{itemize}
\end{proposition}

\begin{proof}
Let us start with item (i). To perform calculations rigorously, it is best to
use the environment representation (\ref{E:rep_int}). Recall also that
$\eta_{k}$ is given by (\ref{D:eta1}). Then
\begin{align}
v_{\ell}  &  =\frac{2}{t}\mbox{{\bf E}}\biggl[\int_{\frac{t}{2}}^{t}%
\int_{\mathbb{R}}\exp(iub_{s})\nu(ds,du)\int_{I_{\ell}}\int_{\frac{t}{2}}%
^{t}\int_{\mathbb{R}}\exp(iux)\nu(ds,du)dx\biggr]\nonumber\\
&  =\frac{2}{t}\int_{I_{\ell}}\mbox{{\bf E}}\biggl[\int_{\frac{t}{2}}^{t}%
\int_{\mathbb{R}}\exp(iub_{s})\nu(ds,du)\int_{\frac{t}{2}}^{t}\int
_{\mathbb{R}}\exp(iux)\nu(ds,du)\biggr]dx.\nonumber
\end{align}
Thanks to (\ref{E:cov_X}), and according to the fact that $\hat{Q}$ is the
Fourier transform of $Q$, we thus have
\begin{align}
v_{\ell}  &  =\frac{2}{t}\int_{I_{\ell}}\left[  \int_{\frac{t}{2}}^{t}%
\int_{\mathbb{R}}\exp(iu(b_{s}-x))\hat{Q}(du)ds\right]  dx\nonumber\\
&  =\frac{2}{t}\int_{\frac{t}{2}}^{t}\int_{I_{\ell}}Q(b_{s}-x)dxds
\label{E:lk}\\
&  \leq\sup_{s\in\lbrack t/2,t]}\int_{I_{\ell}}Q(b_{s}-x)dx. \label{E:ll}%
\end{align}
However, if $\ell\neq k$, on the event $L_{k}$, it is easily checked that, for
$s\in\lbrack t/2,t]$, and for all $x\in I_{\ell}$, we have
\[
\left(  2|\ell-k|-2\right)  t^{\alpha}\leq\left\vert b_{s}-x\right\vert .
\]
According to the fact that $Q$ is a positive decreasing function on
${\mathbb{R}}_{+}$, and $Q\left(  x\right)  =Q\left(  \left\vert x\right\vert
\right)  $, for each $s\in\lbrack t/2,t]$ we can conclude that
\[
\int_{I_{\ell}}Q(b_{s}-x)dx=\int_{I_{\ell}}Q(\left\vert b_{s}-x\right\vert
)dx\leq\int_{t^{\alpha}(2\ell+1)}^{t^{\alpha}(2\ell+1)}Q(\left(
2|\ell-k|-2\right)  t^{\alpha})dx\newline\leq2t^{\alpha}Q(t^{\alpha}%
(2|\ell-k|-2)).\label{E:ant}%
\]
Consequently, putting together equations~(\ref{E:ll}) and~(\ref{E:ant}), we
get
\begin{align}
\left\Vert v\right\Vert _{\tau,k}  &  =v_{k}+\sum_{\ell\neq k}|\ell-k|^{\tau
}v_{\ell}\leq v_{k}+2t^{\alpha}\sum_{\ell\neq k}|\ell-k|^{\tau}Q(t^{\alpha
}(2|\ell-k|-2))\nonumber\\
&  \leq v_{k}+\frac{\kappa}{t^{\alpha(2+\theta)}}\sum_{\ell\neq k}%
|\ell-k|^{-(3+\theta-\tau)}\leq v_{k}+\frac{\kappa}{t^{\alpha(2+\theta)}},
\label{E:tri}%
\end{align}
where $\kappa$ is a positive constant that can change from from one occurence
to the next, and where we have used again Hypothesis \ref{env1}. It is now
readily checked that $\Vert v\Vert_{\tau,k}\leq v_{k}+O(t^{-\alpha})$, which
ends the proof of item (i).\vspace*{0.1in}

Let us prove now item (ii): go back to equation (\ref{E:lk}) and set $\ell=k$.
Then we get
\[
\inf_{s\in\lbrack\frac{t}{2},t]}\int_{I_{k}}Q(b_{s}-x)dx\leq v_{k}\leq
\sup_{s\in\lbrack\frac{t}{2},t]}\int_{I_{k}}Q(b_{s}-x)dx\leq\int_{{\mathbb{R}%
}}Q(u)du=1
\]
To find a lower bound on the left-hand side, we now make use of the
non-degeneracy assumption, as noted in Remark \ref{R:nondeg}: since $Q$ is an
even function, we get $\int_{0}^{\infty}Q\left(  x\right)  dx=1/2$. But if
$b\in L_{k}$, then for any $s\in\lbrack t/2,t]$, we have that the interval
$b_{s}-I_{k}$ contains either $[0,t^{\alpha}]$ or $[-t^{\alpha},0]$, so that,
again by the evenness of $Q$,
\[
\int_{I_{k}}Q(b_{s}-x)dx\geq\int_{0}^{t^{\alpha}}Q\left(  x\right)  dx.
\]
The latter quantity, which tends to $1/2$ when $t\rightarrow\infty$, can be
made to exceed $1/4$ for $t$ large enough. This finishes the proof of item
(ii) with $\underline{c}=1/4$ and $\overline{c}=1$, and the proposition.
\end{proof}

\subsection{Inversion of $C(t)$}

In this section, we will be concerned with the operator $C^{-1}(t)$, where
$C(t)$ has been defined by relation (\ref{covet}), and more specifically, we
will get some information about the solution $\delta$ to the system $C(t)x=v$.
The importance of $\delta$ stems from the fact that the variables $\eta_{k}$
will be independent of $-H_{t}(b)-\sum_{j\in\mathbb{Z}}\delta_{j}\eta_{j}$,
which will be useful for further computations (see Proposition \ref{P:troi}).
However, we have already seen that $C(t)$ behaves asymptotically like the
identity matrix, and thus the vector $\delta$ should be of the same kind as
$v$, in particular when $b\in L_{k}$. This is indeed the case, and will be
proved in the following proposition.

\begin{proposition}
\label{P:nor_del}Under Hypothesis \ref{env1}, suppose in addition that $b\in
L_{k}$. Set $l_{\tau,k}=\{x\in{\mathbb{R}}^{{\mathbb{Z}}};\Vert x\Vert
_{\tau,k}<\infty\}$. Then

\begin{itemize}
\item[(i)] The operator $C(t)$ is invertible in $l_{\tau,k}$. We set then
$\delta:=C^{-1}(t)v$.

\item[(ii)] There exist some strictly positive real numbers $\underline{d}$
and $\overline{d}$ such that
\begin{align}
\underline{d}\leq\delta_{k} \leq\overline{d}.\nonumber
\end{align}

\item[(iii)] The following relation holds:
\[
\left\Vert \delta\right\Vert _{\tau,k}-\delta_{k}=O\left(  \frac{1}{t^{\alpha
}}\right)  \cdot
\]

\item[(iv)] On the probability space $(\Omega,\mathcal{G},\mbox{{\bf P}})$,
the family $\{\eta_{l};l\in{\mathbb{Z}}\}$ is independent of $-H_{t}%
(b)-\sum_{j\in\mathbb{Z}}\delta_{j}\eta_{j}$.
\end{itemize}
\end{proposition}

\begin{remark}
Notice that Proposition \ref{P:nor_del} contains a considerable amount of the
information which will be used for the proof of Lemma \ref{L:env}. Indeed,
inequality (\ref{E:fin}) will be obtained thanks to item (iv), item (iii) will
be invoked for inequality (\ref{whenQcondchanges}), and item (ii) will be
essential in order to define the random variables $\check{\eta}_{0}$ and
$\check{\eta}_{k}$ in (\ref{defchheta}).
\end{remark}

\begin{proof}
[Proof of Proposition \ref{P:nor_del}]~\vspace*{0.1in}

\noindent\emph{Step 1: proving item (i).} We choose the standard operator norm
on $l_{\tau,k}$: a matrix $A$ is defined to be in the linear operator space
$\mathcal{L}_{\tau,k}$ if the norm%
\[
\Vert A\Vert_{\tau,k}:=\sup_{x\in l_{\tau,k}:\left\Vert x\right\Vert _{\tau
,k}=1}\Vert Ax\Vert_{\tau,k}%
\]
is finite. Then, on one hand, the following relations are satisfied since we
are dealing with the operator norm on $l_{\tau,k}$: for $D_{1},D_{2}%
\in\mathcal{L}_{\tau,k}$ and $x\in l_{\tau,k}$:
\begin{equation}
\Vert D_{1}x\Vert_{\tau,k}\leq\Vert D_{1}\Vert_{\tau,k}\Vert x\Vert_{\tau
,k},\quad\mbox{ and }\quad\Vert D_{1}+D_{2}\Vert_{\tau,k}\leq\Vert D_{1}%
\Vert_{\tau,k}+\Vert D_{2}\Vert_{\tau,k}. \label{optk}%
\end{equation}

On the other hand, let us now prove that, setting $A(t):=\mbox{Id}-\lambda
C(t)$, Proposition \ref{P:ord_t} yields that $\Vert A(t)\Vert_{\tau
,k}=O(t^{-\alpha})$, and thus
\begin{equation}
\Vert A(t)\Vert_{\tau,k}<1 \label{Acontraction}%
\end{equation}
if $t$ is large enough. First recall that by definition of $C\left(  t\right)
$ and $\lambda$, denoting by $\dot{C}\left(  t\right)  $ the matrix $C\left(
t\right)  $ deprived of its diagonal, we have%
\[
A\left(  t\right)  =-\lambda\dot{C}\left(  t\right)  .
\]
By Proposition \ref{P:ord_t} item (i), $\lambda$ tends to $1$ as
$t\rightarrow\infty$. Therefore, it is sufficient to show that $\Vert\dot
{C}(t)\Vert_{\tau,k}=O(t^{-\alpha})$. Thus let $x\in l_{\tau,k}$ such that
$\left\Vert x\right\Vert _{\tau,k}=1$. In other words,%
\[
\left\vert x_{k}\right\vert +\sum_{i\neq k}\left\vert x_{i}\right\vert
\left\vert i-k\right\vert ^{\tau}=1.
\]
Now we calculate the two terms that form $\Vert\dot{C}(t)x\Vert_{\tau,k}$. The
first is%
\begin{align}
\left\vert \left(  \dot{C}(t)x\right)  _{k}\right\vert  &  =\left\vert
\sum_{j\neq k}C_{kj}\left(  t\right)  x_{j}\right\vert \leq\sum_{j\neq
k}\left\vert C_{kj}\left(  t\right)  x_{j}\right\vert \nonumber\\
&  \leq\left(  \sum_{j\neq k}\left\vert x_{j}\right\vert \left\vert
k-j\right\vert ^{\tau}\right)  \left(  \sum_{j\neq k}\left\vert C_{kj}\left(
t\right)  \right\vert \left\vert k-j\right\vert ^{\tau}\right) \nonumber\\
&  \leq1\cdot O\left(  t^{-\alpha}\right)  \label{neumannorm}%
\end{align}
where we used the assumption $\left\Vert x\right\Vert _{\tau,k}=1$ and the
result of Proposition \ref{P:ord_t} item (ii). The second term in $\Vert
\dot{C}(t)x\Vert_{\tau,k}$ equals
\[
\sum_{i\neq k}\left\vert \sum_{j\neq i}C_{ij}\left(  t\right)  x_{j}%
\right\vert \left\vert i-k\right\vert ^{\tau}\leq\sum_{j\in{\mathbb{Z}}%
}\left\vert x_{j}\right\vert \sum_{i\neq j;i\neq k}\left\vert C_{ij}\left(
t\right)  \right\vert \left\vert i-k\right\vert ^{\tau}=:K_{2};
\]
we split this sum up according to $j=k$ or $j\neq k$:%
\begin{align*}
K_{2}  &  \leq\left\vert x_{k}\right\vert \sum_{i\neq k}\left\vert
C_{ik}\left(  t\right)  \right\vert \left\vert i-k\right\vert ^{\tau}%
+\sum_{j\neq k}\left\vert x_{j}\right\vert \sum_{i\neq j;i\neq k}\left\vert
C_{ij}\left(  t\right)  \right\vert \left\vert i-j+j-k\right\vert ^{\tau}\\
&  \leq\left\vert x_{k}\right\vert \sum_{i\neq k}\left\vert C_{ik}\left(
t\right)  \right\vert \left\vert i-k\right\vert ^{\tau}+\sum_{j\neq
k}\left\vert x_{j}\right\vert \sum_{i\neq j;i\neq k}\left\vert C_{ij}\left(
t\right)  \right\vert \left\vert i-j\right\vert ^{\tau}\\
&  +\sum_{j\neq k}\left\vert x_{j}\right\vert \left\vert j-k\right\vert
^{\tau}\sum_{i\neq j;i\neq k}\left\vert C_{ij}\left(  t\right)  \right\vert
\end{align*}
where in the last line we used the fact that $|a+b|^{\tau}\leq|a|^{\tau
}+|b|^{\tau}$ whenever $\tau\in(0,1)$.

Now using the fact that $\sum_{i\neq j;i\neq k}\left\vert C_{ij}\left(
t\right)  \right\vert $ is bounded above by $\sum_{i\neq j}\left\vert
C_{ij}\left(  t\right)  \right\vert \left\vert i-j\right\vert ^{\tau}$, and
the latter is $O\left(  t^{-\alpha}\right)  $ by Proposition \ref{P:ord_t}
item (ii), we can assert $K_{2}\leq O\left(  t^{-\alpha}\right)  $, which,
combined with (\ref{neumannorm}), implies our goal $\Vert\dot{C}(t)\Vert
_{\tau,k}=O(t^{-\alpha})$, and thus (\ref{Acontraction}). This contraction
relation (\ref{Acontraction}) finishes the proof of (i) because it allows us
to define $C^{-1}(t)$ in $\mathcal{L}_{\tau,k}$ by a Von Neumann type series
of the form
\begin{equation}
C^{-1}(t)=\lambda\sum_{j\geq0}A^{j}. \label{E:VN}%
\end{equation}
\vspace*{0.1in}

\noindent\emph{Step 2: proving item (ii).} For $t$ large enough, set
$\delta=C^{-1}(t)v$, which makes sense since $v\in l_{\tau,k}$. Then, thanks
to the fact that $C^{-1}(t)$ can be defined by relation (\ref{E:VN}), we have
\begin{align*}
\delta_{k}=\lambda\Bigl(v_{k}+\sum_{j\geq1}(A^{j}v)_{k}\Bigr) &  \geq
\lambda\Bigl(v_{k}-\sum_{j\geq1}\Vert A^{j}v\Vert_{\tau,k}\Bigr)\\
&  \geq\lambda\Bigl(v_{k}-\sum_{j\geq1}\Vert A\Vert_{\tau,k}^{j}\Vert
v\Vert_{\tau,k}\Bigr),
\end{align*}
where we have used the relations $x_{k}\geq-\Vert x\Vert_{\tau,k}$ and
(\ref{optk}). Hence, since $\Vert A(t)\Vert_{\tau,k}=O(t^{-\alpha})$, we
obtain
\begin{equation}
\delta_{k}\geq\lambda\left(  v_{k}-\frac{\Vert A\Vert_{\tau,k}}{1-\Vert
A\Vert_{\tau,k}}\Vert v\Vert_{\tau,k}\right)  \geq\lambda\left(
v_{k}+O\left(  \frac{1}{t^{\alpha}}\right)  \right)  \geq\underline
{d}+O\left(  \frac{1}{t^{\alpha}}\right)  ,\label{E:di}%
\end{equation}
according to the properties of $v$ shown at Proposition \ref{P:nor_v}. The
upper bound on $\delta_{k}$ can now be shown by the same type of argument,
which ends the proof of our claim. \vspace*{0.1in}

\noindent\emph{Step 3: proving item (iii).} Let us evaluate now the quantity
$\Vert\delta\Vert_{\tau,k}-\delta_{k}$: thanks to relations (\ref{optk}) and
(\ref{E:di}), we get
\begin{align}
\Vert\delta\Vert_{\tau,k}-\delta_{k}  &  \leq\Vert C(t)^{-1}\Vert_{\tau
,k}\Vert v\Vert_{\tau,k}-\delta_{k}\nonumber\\
&  \leq\left(  \Vert C(t)^{-1}\Vert_{\tau,k}\Vert v\Vert_{\tau,k}-\lambda
v_{k}+\frac{\lambda\Vert A\Vert_{\tau,k}}{1-\Vert A\Vert_{\tau,k}}\Vert
v\Vert_{\tau,k}\right)  .\nonumber
\end{align}
Thus, using again that fact that $C^{-1}(t)$ is defined by
equation~(\ref{E:VN}) and relation (\ref{optk}), we obtain
\begin{align}
\Vert\delta\Vert_{\tau,k}-\delta_{k}  &  \leq\lambda\left(  \frac{1+\Vert
A\Vert_{\tau,k}}{1-\Vert A\Vert_{\tau,k}}\Vert v\Vert_{\tau,k}-v_{k}\right)
\nonumber\\
&  =\lambda\left(  \Vert v\Vert_{\tau,k}-v_{k}\right)  +O\left(  \frac
{1}{t^{\alpha}}\right)  =O\left(  \frac{1}{t^{\alpha}}\right)  ,\nonumber
\end{align}
where in the last two steps, we have invoked, respectively, item (i) and
Proposition~\ref{P:nor_v}. This concludes our proof of (iii).\vspace*{0.1in}

\noindent\emph{Step 4: proving item (iv).} Recall that, by definition,
$C(t)=t^{-\left(  1-\alpha\right)  }\mbox{\bf Cov}(\eta)$. Hence
\begin{align*}
\delta_{j}  &  =\left(  C^{-1}\left(  t\right)  v\right)  _{j}=\frac{1}%
{4}t^{1-\alpha}\sum_{k\in\mathbb{Z}}\left[  \mbox{\bf Cov}(\eta)\right]
_{jk}^{-1}v_{k}\\
&  =\sum_{k\in\mathbb{Z}}\left[  \mbox{\bf Cov}(\eta)\right]  _{jk}%
^{-1}\mathbf{E}\left[  \int_{\frac{t}{2}}^{t}W(ds,b_{s})\,\eta_{k}\right] \\
&  =\sum_{k\in\mathbb{Z}}\left[  \mbox{\bf Cov}(\eta)\right]  _{jk}%
^{-1}\mathbf{E}\left[  \left(  -H_{t}\left(  b\right)  \right)  \eta
_{k}\right]  ,
\end{align*}
\text{we have the following standard calculation for any }$\ell\in\mathbb{Z}$%
\begin{align*}
&  \mathbf{E}\left[  \left(  -H_{t}(b)-\sum_{j\in\mathbb{Z}}\delta_{j}\eta
_{j}\right)  \eta_{\ell}\right] \\
&  =-\mathbf{E}\left[  H_{t}\left(  b\right)  \eta_{\ell}\right]
+\mathbf{E}\sum_{j\in\mathbb{Z}}\sum_{k\in\mathbb{Z}}\left[
\mbox{\bf Cov}(\eta)\right]  _{jk}^{-1}\mathbf{E}\left[  H_{t}\left(
b\right)  \eta_{k}\right]  \eta_{j}\eta_{\ell}\\
&  =-\mathbf{E}\left[  H_{t}\left(  b\right)  \eta_{\ell}\right]
+\mathbf{E}\sum_{j\in\mathbb{Z}}\sum_{k\in\mathbb{Z}}\left[
\mbox{\bf Cov}(\eta)\right]  _{jk}^{-1}\left[  \mbox{\bf Cov}(\eta)\right]
_{j\ell}\mathbf{E}\left[  H_{t}\left(  b\right)  \eta_{k}\right] \\
&  =-\mathbf{E}\left[  H_{t}\left(  b\right)  \eta_{\ell}\right]  +\sum
_{k\in\mathbb{Z}}\delta_{k\ell}\mathbf{E}\left[  H_{t}\left(  b\right)
\eta_{k}\right]  =0.
\end{align*}
Now since for fixed $b$, $H_{t}\left(  b\right)  $ and the sequence $\eta$ are
both linear functionals of a same Gaussian field, they form a jointly Gaussian
vector, and are thus independent.
\end{proof}

\section{Application of Girsanov's theorem\label{girsanov}}

In our context, the cost of having $b$ living in the interval $I_{k}%
=[t^{\alpha}(2k-1),t^{\alpha}(2k+1)]$ instead of $I_{0}=[-t^{\alpha}%
,t^{\alpha}]$ can be calculated explicitly thanks to Girsanov's theorem: given
an integer $k$, a real number $t$ and a realization of the environment $W$, we
define a new environment by setting $W^{k,t}(ds,x):=W(ds,x+h\left(  s\right)
)$, where%
\[
h\left(  s\right)  :=\min(2s/t,1)2kt^{\alpha},
\]
or more rigorously,%
\begin{equation}
W^{k,t}(s,x):=\int_{0}^{s}W(du,x+h\left(  u\right)  ). \label{shiftenv}%
\end{equation}

A simple and useful result that we can now prove is the following.

\begin{lemma}
\label{WWkt}The random fields defined by $W=\{W\left(  s,x\right)  :\left(
s,x\right)  \in{\mathbb{R}}_{+}\times{\mathbb{R\}}}$ and $W^{k,t}=\{\int
_{0}^{s}W(du,u+h\left(  u\right)  ):(s,x)\in{\mathbb{R}}_{+}\times
{\mathbb{R\}}}$ have the same distribution.
\end{lemma}

\begin{proof}
The easiest way to establish this result is to revert to the representation of
$W$ using the Gaussian measure $\nu$, i.e. (\ref{E:rep_int}), and also its
consequence (\ref{E:rep_intH}), so that
\[
W^{k,t}(s,x):=\int_{0}^{s}\int_{{\mathbb{R}}}e^{\imath\lambda(x+h(u))}%
\nu(ds,du).
\]
Since the law of this centered Gaussian field is determined by its covariance
structure only, it is now immediate to check, using the formulas
(\ref{defnuf}) and (\ref{E:cov_X}), that it has the same law as $W$, since we
have%
\[
W(s,x):=\int_{0}^{s}\int_{{\mathbb{R}}}e^{\imath\lambda x}\nu(ds,du).
\]
The calculations are left to the reader.
\end{proof}

\begin{proof}
[Alternate Proof]It is also possible to invoke a direct proof of this fact,
using $L^{2}$ approximations of $W^{k,t}(s,x)$ by Riemann sums. For fixed
$s,x$, $W^{k,t}\left(  s,x\right)  $ can be written as a limit in
$L^{2}\left(  \Omega\right)  $, as $n\rightarrow\infty$, of the sum
$\sum_{i=1}^{n}J_{i}^{k,t}$ of the increments $J_{i}^{k,t}%
:=W([si/n,s(i+1)/n],x+h\left(  si/n\right)  )$, whose individual laws are
identical to those of the $J_{i}$'s defined without adding the shift $h\left(
si/n\right)  $, because $W$ is spatially homogeneous. Since the $J_{i}^{k,t}%
$'s are independent as $i$ changes, (as are the $J_{i}$'s), $W^{k,t}\left(
s,x\right)  $ and $W\left(  s,x\right)  $ have the same distribution for fixed
$s,x$; we omit the end of this -- more intuitive but less rigorous -- proof.
\end{proof}

We also need to introduce a modified partition function $\tilde{Z}$ defined
by
\begin{equation}
\tilde{Z}_{t}^{\alpha}(k)=E_{b}\biggl[\mathbf{1}_{L_{k}}\left(  b\right)
\exp\biggl(\beta\biggl(\int_{0}^{t}W(ds,b_{s})-\sum_{j\in\mathbb{Z}}\delta
_{j}\eta_{j}\biggr)\biggr)\biggr]. \label{Ztildadef}%
\end{equation}
In the sequel, we will have to stress the dependence of these partition
functions on the environment under consideration. We will thus set $\tilde
{Z}_{t}^{\alpha}(k)=\tilde{Z}_{t}^{\alpha}(k,W)$. With these notations in
mind, we can prove the following proposition, which shows that the cost of
having $b$ live in $L_{k}$ rather than $L_{0}$ is exponential of order
$t^{2\alpha-1}$.

\begin{proposition}
\label{P:entr} Given two positive real numbers $\alpha$ and $t$, and an
integer $k$ fixed, we have%
\begin{equation}
\tilde{Z}_{t}^{\alpha}(k,W)\geq\exp\left[  -4(k+k^{2})t^{2\alpha-1}\right]
~\tilde{Z}_{t}^{\alpha}(0,W^{k,t}). \label{corgirs}%
\end{equation}

\end{proposition}

\begin{proof}
~\vspace*{0.1in}

\noindent\emph{Step 1: using Girsanov's theorem.} Given $k$ and $t$, and with
$h\left(  s\right)  =\min(2s/t,1)2kt^{\alpha}$ as defined above, we associate
to a path $b$ a shifted path $b^{\prime}$ by the relation
\[
b_{s}^{\prime}\equiv b_{s}-h\left(  s\right)  ,\quad\mbox{ for }\quad
s\in{\mathbb{R}}.
\]
Notice that this shift transforms a path which lives in the interval $I_{k}$
for all $s\in\lbrack t/2,t]$ into a path which belongs to $I_{0}$ in the same
time interval. More precisely, one immediately checks that $\mathbf{1}_{L_{k}%
}\left(  b\right)  =\mathbf{1}_{L_{0}}\left(  b^{\prime}\right)  $. Let us
call $M_{t}(b^{\prime})$ the Girsanov density involved in the shift between
$b$ and $b^{\prime}$, that is
\[
M_{t}(b^{\prime})=\exp\left(  -b_{t/2}^{\prime}4kt^{\alpha-1}-4k^{2}%
t^{2\alpha-1}\right)  .
\]
The choice of $h\left(  s\right)  =4kst^{\alpha-1}$ for $s\in\lbrack0,t/2]$ is
made to obtain a continuous function that starts at $0$, and is piecewise
linear (constant over $[t/2,t]$); this function has the advantage that its
Girsanov \textquotedblleft energy\textquotedblright\ is minimal, ensuring that
our proof is most efficient. It is possible that other, non-linear, choices
could have fulfilled our purposes, but this would be an unnecessary complication.

For sake of clarity, let us stress now the dependence of the random variables
$\delta,\eta$, etc., on the data of our problem: it is readily checked for
instance that
\[
\eta_{j}=\eta_{j}\left(  W\right)  ,\quad\mbox{ and }\quad\delta_{j}%
=\delta_{j}\left(  b,\mathcal{L}(W)\right)  ,
\]
where a function of $\left(  W\right)  $ represents its dependence on the
increments of $W$ in the interval $[0,t]$, as a random variable, where the
symbol $\mathcal{L}\left(  \cdot\right)  $ denotes the law (distribution) of a
process on $[0,t]$, and where a function of $b$ represents its dependence on
the fixed path $b$. Then, adopting this convention, we have
\[
\tilde{Z}_{t}^{\alpha}(k,W)=E_{b}\Bigg[\mathbf{1}_{L_{k}}\exp\Bigg(\beta
\int_{0}^{t}W(ds,b_{s}^{\prime}+h\left(  s\right)  )-\sum_{j\in{\mathbb{Z}}%
}\delta_{j}\left(  b^{\prime}+h,\mathcal{L}(W)\right)  \eta_{j}%
(W)\Bigg)\Bigg].
\]
After applying Girsanov's transformation, noting that by definition, $\int
_{0}^{t}W(ds,b_{s}^{\prime}+h\left(  s\right)  )=\int_{0}^{t}W^{k,t}%
(ds,b_{s}^{\prime})$, we get (recall that $b^{\prime}$ is a standard Brownian
motion under the new probability, so that it is notationally legitimate to
write $b$ instead of $b^{\prime}$, and to denote expectation with respect to
the new measure by $E_{b}$):
\[
\tilde{Z}_{t}^{\alpha}(k,W)=E_{b}\bigg[\mathbf{1}_{L_{0}(b)}M_{t}%
(b)\exp\biggl(\beta\biggl(\int_{0}^{t}W^{k,t}(ds,b_{s})r-\sum_{j\in\mathbb{Z}%
}\delta_{j}\left(  b+h,\mathcal{L}(W)\right)  \eta_{j}%
(W)\biggr)\biggr)\bigg].\label{eq:Zt-alpha-girsanov}%
\]
\vspace*{0.1in}

\noindent\emph{Step 2: reexpressing the transformed }$\eta$. One should now
compare the random variables $\eta_{j}(W)$ and $\eta_{j}(W^{k,t})$: by
definition of these quantities, we have
\begin{align}
\eta_{j}(W^{k,t}) &  =\frac{1}{t^{2\alpha}}\int_{(2j-1)t^{\alpha}%
}^{(2j+1)t^{\alpha}}\int_{t/2}^{t}W(ds,x+2kt^{\alpha})dx\nonumber\\
&  =\frac{1}{t^{2\alpha}}\int_{(2(j+k)-1)t^{\alpha}}^{(2(j+k)+1)t^{\alpha}%
}\int_{t/2}^{t}W(ds,x)dx=\eta_{j+k}(W).\label{eq:eta-j-shifted}%
\end{align}
In particular, the law of $\eta(W^{k,t})$, considered as the set of random
variables forming that sequence, is the same as the law of $\eta(W)$, a fact
which we will not use in this proof, but will be crucial in the proof of the
next lemma. \vspace*{0.1in}

\noindent\emph{Step 3: reexpressing the transformed }$\delta$. Along the same
lines as (\ref{eq:eta-j-shifted}), we now show that
\begin{equation}
\delta_{j}\left(  b+h,\mathcal{L}(W)\right)  =\delta_{j-k}\left(
b,\mathcal{L}(W^{k,t})\right)  .\label{eq:delta-j-shifted}%
\end{equation}
To see this, we recall the definition of $\delta$: we have
\[
\delta=\delta\left(  b+h,\mathcal{L}\left(  W\right)  \right)  =\left[
C\left(  t\right)  \right]  ^{-1}v=\left[  C\left(  t,\mathcal{L}\left(
W\right)  \right)  \right]  ^{-1}v\left(  b+h,\mathcal{L}\left(  W\right)
\right)  ,
\]
where we calculate
\begin{align*}
&  C_{\ell,m}\left(  t,\mathcal{L}\left(  W\right)  \right)  \\
&  =\frac{1}{t^{(\alpha+1)}}\mathbf{E}\left[  \int_{\frac{t}{2}}^{t}%
\int_{(2m-1)t^{\alpha}}^{(2m+1)t^{\alpha}}W(ds,x)dx\cdot\int_{\frac{t}{2}}%
^{t}\int_{(2\ell-1)t^{\alpha}}^{(2\ell+1)t^{\alpha}}W(ds,x)dx\right]  \\
&  =\frac{1}{t^{(\alpha+1)}}\mathbf{E}\left[  \int_{\frac{t}{2}}^{t}%
\int_{(2(m-k)-1)t^{\alpha}}^{(2(m-k)+1)t^{\alpha}}W(ds,x+2kt^{\alpha}%
)dx\cdot\int_{\frac{t}{2}}^{t}\int_{(2(\ell-k)-1)t^{\alpha}}^{(2(\ell
-k)+1)t^{\alpha}}W(ds,x+2kt^{\alpha})dx\right]  \\
&  =\frac{1}{t^{(\alpha+1)}}\mathbf{E}\left[  \int_{\frac{t}{2}}^{t}%
\int_{I_{m-k}}W^{k,t}(ds,x)dx\cdot\int_{\frac{t}{2}}^{t}\int_{I_{\ell-k}%
}W^{k,t}(ds,x)dx\right]  \\
&  =C_{\ell-k,m-k}\left(  t,\mathcal{L}\left(  W^{k,t}\right)  \right)  ,
\end{align*}
and similarly%
\begin{align*}
v_{\ell}\left(  b+h,\mathcal{L}\left(  W\right)  \right)   &  =4t^{\alpha
-1}\mathbf{E}\left[  \int_{\frac{t}{2}}^{t}\int_{(2\ell-1)t^{\alpha}}%
^{(2\ell+1)t^{\alpha}}W(ds,x)dx\cdot\int_{\frac{t}{2}}^{t}W(ds,b_{s}+h\left(
s\right)  )\right]  \\
&  =4t^{\alpha-1}\mathbf{E}\left[  \int_{\frac{t}{2}}^{t}\int_{(2(\ell
-k)-1)t^{\alpha}}^{(2(\ell-k)+1)t^{\alpha}}W(ds,x+h\left(  s\right)
)dx\cdot\int_{\frac{t}{2}}^{t}W(ds,b_{s}+h\left(  s\right)  )\right]  \\
&  =4t^{\alpha-1}\mathbf{E}\left[  \int_{I_{\ell}}\int_{\frac{t}{2}}%
^{t}W^{k,t}(ds,x)dx\cdot\int_{\frac{t}{2}}^{t}W^{k,t}(ds,b_{s})\right]  \\
&  =v_{\ell-k}\left(  b,\mathcal{L}\left(  W^{k,t}\right)  \right)  .
\end{align*}
We may thus write that the definition of $\delta\left(  b+h,\mathcal{L}\left(
W\right)  \right)  $ is equivalent to,%
\begin{align*}
\forall\ell &  \in{\mathbb{Z}}:\sum_{m\in{\mathbb{Z}}}C_{\ell,m}\left(
t,\mathcal{L}\left(  W\right)  \right)  \delta_{m}\left(  b+h,\mathcal{L}%
\left(  W\right)  \right)  =v_{\ell}\left(  b+h,\mathcal{L}\left(  W\right)
\right)  \\
&  \iff\forall\ell\in{\mathbb{Z}}:\sum_{m\in{\mathbb{Z}}}C_{\ell-k,m-k}\left(
t,\mathcal{L}\left(  W^{k,t}\right)  \right)  \delta_{m}\left(
b+h,\mathcal{L}\left(  W\right)  \right)  =v_{\ell-k}\left(  b,\mathcal{L}%
\left(  W^{k,t}\right)  \right)  \\
&  \iff\forall\ell\in{\mathbb{Z}}:\sum_{m\in{\mathbb{Z}}}C_{\ell,m}\left(
t,\mathcal{L}\left(  W^{k,t}\right)  \right)  \delta_{m+k}\left(
b+h,\mathcal{L}\left(  W\right)  \right)  =v_{\ell}\left(  b,\mathcal{L}%
\left(  W^{k,t}\right)  \right)  .
\end{align*}
This last statement is equivalent to saying $\delta_{m+k}\left(
b+h,\mathcal{L}\left(  W\right)  \right)  =\delta_{m}\left(  b,\mathcal{L}%
\left(  W^{k,t}\right)  \right)  $, which is precisely the statement of
(\ref{eq:delta-j-shifted}).\vspace*{0.1in}

\noindent\emph{Step 4: conclusion.} Plugging equations (\ref{eq:eta-j-shifted}%
) and (\ref{eq:delta-j-shifted}) into (\ref{eq:Zt-alpha-girsanov}), we end up
with
\[
\tilde{Z}_{t}^{\alpha}(k,W)=E_{b}\bigg[\mathbf{1}_{L_{0}(b)}M_{t}%
(b)\exp\biggl(\beta\biggl(\int_{0}^{t}W^{k,t}(ds,b_{s})-\sum_{j\in\mathbb{Z}%
}\delta_{j-k}\left(  b,\mathcal{L}(W^{k,t})\right)  \eta_{j-k}(W^{k,t}%
)\biggr)\biggr)\bigg].\label{girsanovproofconclu}%
\]

To conclude the proof of the proposition, notice that for $b\in L_{0}$, we get
$\left\vert b_{t/2}\right\vert \leq t^{\alpha}$, and therefore%
\begin{equation}
M_{t}\left(  b\right)  \geq\exp\left(  -4kt^{2\alpha-1}-4k^{2}t^{2\alpha
-1}\right)  . \label{girsanovproofconclu2}%
\end{equation}
Combining (\ref{girsanovproofconclu}) and (\ref{girsanovproofconclu2}), and
renumbering the sum for $j\in\mathbb{Z}$ as $j^{\prime}=j-k\in\mathbb{Z}$, we
recognize the term $\tilde{Z}_{t}^{\alpha}(0,W^{k,t})$, and the proof is complete.
\end{proof}

The above proof has an important consequence which we record here for use at a
crucial point in the next section.

\begin{lemma}
\label{independence}Let
\[
X\left(  W,b\right)  =-H_{t}(b)-\sum_{j\in\mathbb{Z}}\delta_{j}\eta_{j}%
=\int_{0}^{t}W(ds,b_{s})-\sum_{j\in\mathbb{Z}}\delta_{j}\left(  b,\mathcal{L}%
(W)\right)  \eta_{j}\left(  W\right)
\]
and therefore
\[
X\left(  W^{k,t},b\right)  =\int_{0}^{t}W^{k,t}(ds,b_{s})-\sum_{j\in
\mathbb{Z}}\delta_{j}\left(  b,\mathcal{L}(W^{k,t})\right)  \eta_{j}\left(
W^{k,t}\right)  .
\]
Denote by $\eta\left(  W\right)  $ the entire sequence $\left\{  \eta
_{j}\left(  W\right)  :j\in{\mathbb{Z}}\right\}  $. Then for each $b$,
$X\left(  W,b\right)  $ and $\eta\left(  W\right)  $ are independent, and for
each $k\in{\mathbb{Z}}$, and each $b$, $X\left(  W^{k,t},b\right)  $ and
$\eta\left(  W\right)  $ are independent.
\end{lemma}

\begin{proof}
We have already proved in Proposition \ref{P:nor_del} (iv) that $X\left(
W,b\right)  $ and $\eta\left(  W\right)  $ are independent, which is the first
half of what we have to prove. This implies in addition that $X\left(
W^{k,t},b\right)  $ and $\eta\left(  W^{k,t}\right)  $ are also independent
because the random fields $W\ $and $W^{k,t}$ have the same distribution (Lemma
\ref{WWkt}).

To conclude the proof this lemma, we simply invoke the portion of the proof of
Proposition \ref{P:entr} which shows the specific shift equality relation
$\eta_{j+k}\left(  W\right)  =\eta_{j}\left(  W^{k,t}\right)  $, from
(\ref{eq:eta-j-shifted}): this is a $\mathbf{P}$-almost-sure equality in
$\Omega$. This implies that the sets of points in the sequences $\left\{
\eta_{j}\left(  W\right)  :j\in\mathbb{Z}\right\}  $ and $\left\{  \eta
_{j}\left(  W^{k,t}\right)  :j\in\mathbb{Z}\right\}  $ are precisely the same
sets of random variables. Therefore, for each $k$ and $b$, $X\left(
W^{k,t},b\right)  $ is independent of the entire sequence $\eta\left(
W\right)  $.
\end{proof}

\section{Proof of Lemma \ref{L:env}\label{Lenvproof}}

Recall that we have reduced our problem to the evaluation of
$\mbox{{\bf P}}(\mathcal{B}_{t})$, where
\[
\mathcal{B}_{t}=\mathcal{A}_{t}^{c}=\left\{  \mbox{For all }k\in{\mathbb{Z}%
},\,Z_{t}^{\alpha}(k)\leq Z_{t}^{\alpha}(0)\right\}  ,
\]
and one wishes to show that $\lim_{t\rightarrow\infty}%
\mbox{{\bf P}}(\mathcal{B}_{t})=0$. Then a first step in order to prove this
claim is to truncate $\mathcal{B}_{t}$: for a positive integer $M$ let
$\mathbb{\dot{Z}}_{M}$ and $\bar{\mathbb{Z}}_{M}$ be the sets defined
respectively by
\begin{equation}
\bar{\mathbb{Z}}_{M}=\{-M,-M+1,\ldots,M-1,M\}\quad\mbox{ and }\quad
\mathbb{\dot{Z}}_{M}=\bar{\mathbb{Z}}_{M}\backslash\{0\}, \label{defbzm}%
\end{equation}
and $\mathcal{B}_{M,t}$ the event defined by
\[
\mathcal{B}_{M,t}=\left\{  \mbox{For all }k\in\mathbb{\dot{Z}}_{M}%
,\,Z_{t}^{\alpha}(k)\leq Z_{t}^{\alpha}(0)\right\}  .
\]
Then obviously, $\mbox{{\bf P}}(\mathcal{B}_{t})\leq\mbox{{\bf P}}(\mathcal{B}%
_{M,t})$, and we only need to prove that $\mbox{{\bf P}}(\mathcal{B}_{M,t})$
tends to 0 as $t\rightarrow\infty$.

\vspace{0.3cm}

Here is a brief account on the strategy we will follow in order to complete
our proof.

\vspace{0.3cm}

\noindent\textbf{(1)} Recall that we are trying to bound
\begin{equation}
\mbox{{\bf P}}\left(  \mathcal{B}_{M,t}\right)  =\mbox{{\bf P}}\left(
E_{b}\left[  \mathbf{1}_{L_{k}}e^{-\beta H_{t}(b)}\right]  <E_{b}\left[
\mathbf{1}_{L_{0}}e^{-\beta H_{t}(b)}\right]  \mbox{ for all }k\in
\mathbb{\dot{Z}}_{M}\right)  . \label{abnd}%
\end{equation}
A natural idea is then to split the conditions $E_{b}[\mathbf{1}_{L_{k}%
}e^{-\beta H_{t}(b)}]<E_{b}[\mathbf{1}_{L_{0}}e^{-\beta H_{t}(b)}]$ in terms
of a condition involving the random variables $\eta_{l}$ introduced at
(\ref{D:eta1}), on which we have a reasonable control, and another set of
conditions involving some random variables independent of the family
$\{\eta_{l};l\in{\mathbb{Z}}\}$. However, we have already seen in Proposition
\ref{P:nor_del} that $-H_{t}(b)-\sum_{j\in\mathbb{Z}}\delta_{j}\eta_{j}$ is
independent of $\{\eta_{l};l\in{\mathbb{Z}}\}$. Thus, a natural choice will be
to replace $e^{-\beta H_{t}(b)}$ by $e_{t}(b)$ in the expression (\ref{abnd}),
where $e_{t}(b)$ is defined by
\[
e_{t}(b):=\exp\left(  -\beta\left(  H_{t}(b)+\sum_{j\in\mathbb{Z}}\delta
_{j}\eta_{j}\right)  \right)  .
\]
Of course, this induces a correction term $\exp(\beta\sum_{j\in\mathbb{Z}%
}\delta_{j}\eta_{j})$, but this term can be controlled, since the covariance
structure of the family $\{\eta_{l};l\in{\mathbb{Z}}\}$ is given by
Proposition \ref{covet}, and the vector $\delta$ is controlled by means of
Proposition \ref{P:nor_del}. Up to a negligible term, we will be allowed to
bound $\mbox{{\bf P}}(\mathcal{B}_{M,t})$ by a probability of the form
\begin{equation}
\mbox{{\bf P}}\Biggl(\mbox{For any }k\in\mathbb{\dot{Z}}_{M};\,\frac{\tilde
{Z}_{t}^{\alpha}(k)}{\tilde{Z}_{t}^{\alpha}(0)}<\exp(2\gamma t^{2\alpha
-1}+\eta_{k}^{\ast})\Biggr), \label{apbtm}%
\end{equation}
where $\tilde{Z}_{t}^{\alpha}(k)=E_{b}[\mathbf{1}_{L_{k}}e_{t}(b)]$, as was
defined in Section \ref{girsanov} on Girsanov's theorem, the term
$t^{2\alpha-1}$ comes from the sharp estimates of $\delta$ in Proposition
\ref{P:nor_del}, and the random variable $\eta_{k}^{\ast}$ is one which is
defined using only the random variables $\eta$, because it results from using
$e_{t}(b)$ instead of $e^{-H_{t}(b)}$. The effect of $\eta_{k}^{\ast}$ can be
studied separately from the behavior of the ratio $\tilde{Z}_{t}^{\alpha
}(k)/\tilde{Z}_{t}^{\alpha}(0)$, by the independence property of these two quantities.

\vspace{0.3cm}

\noindent\textbf{(2)} Notice that up to now, we have chosen our parameters
carefully in order to get a penalization of order $\exp(2\gamma t^{2\alpha
-1})$ in (\ref{apbtm}). This was chosen to be consistent with the correction
$\exp(-4(k+k^{2})t^{2\alpha-1})$ we must impose on $b$ if we wish that it live
the second half of its life in $I_{k}$, as we showed by using Girsanov's
theorem in Proposition \ref{P:entr}. In fact, we will be able to bound
$\mbox{{\bf P}}(\mathcal{B}_{M,t})$ by $\mbox{{\bf P}}(F_{M})$, where the
event $F_{M}$ is defined by
\[
F_{M}=\left\{  \mbox{For any }k\in\mathbb{\dot{Z}}_{M};\,\frac{\tilde{Z}%
_{t}^{\alpha}(0,W^{k,t})}{\tilde{Z}_{t}^{\alpha}(0,W)}<\exp(\hat{\gamma
}t^{2\alpha-1}+\eta_{k}^{\ast})\right\}  ,
\]
for some constant $\hat{\gamma}=\hat{\gamma}(M)$, where the shifted
environments $W^{k,t}$ are defined in (\ref{shiftenv}).

\vspace{0.3cm}

\noindent\textbf{(3)} It turns out that the random variable $\eta_{k}^{\ast}$
is optimally chosen to be of the order $\eta_{0}-\eta_{k}$ (see the definition
(\ref{defchheta}) we chose below). We are now considering a set $F_{M}$
involving the random variables $\tilde{Z}_{t}$ and $\eta_{k}^{\ast}$, and this
will allow us to take advantage of the following facts:

\begin{enumerate}
\item The ratio $\tilde{Z}_{t}^{\alpha}(0,W^{k,t})/\tilde{Z}_{t}^{\alpha
}(0,W)$ cannot be too small at many different sites $k\in\mathbb{\dot{Z}}_{M}%
$, by translation invariance in space of $W$.

\item Proposition \ref{covet} asserts that $\{t^{-\left(  1-\alpha\right)
/2}\eta_{k};k\in\mathbb{\dot{Z}}_{M}\}$ is asymptotically a standard Gaussian
vector. Since $\eta_{k}^{\ast}$ is of the order $\eta_{0}-\eta_{k}$ (and thus
of magnitude $t^{\left(  1-\alpha\right)  /2}$), it can be highly negative at
many different sites; thus we are allowed to expect that $\exp(\hat{\gamma
}t^{2\alpha-1}+\eta_{k}^{\ast})$ is much smaller than 1 at many different
sites of $\mathbb{\dot{Z}}_{M}$.

\item The random variables $\tilde{Z}_{t}^{\alpha}$ are independent of
anything defined using $\eta$, inlcuding $\eta_{k}^{\ast}$, and hence the two
effects alluded to above can be taken into account separately.
\end{enumerate}

\vspace{0.3cm}

\noindent\textbf{(4)} These heuristic considerations will be formalized in
Step 3 of the proof below, through the introduction of an intricate family of
subsets of $\mathbb{\dot{Z}}_{M}$, but let us mention that the exponent $3/5$
comes out already at this stage: indeed, the above considerations only make
sense if the magnitude $t^{(1-\alpha)/2}$ of the $\eta_{k}^{\ast}$ is greater
than the magnitude $t^{2\alpha-1}$ of the penalization, so that a highly
negative $\eta_{k}^{\ast}$ can win against the latter. This can only occur,
obviously, whenever $\alpha<3/5$. In this sense, our estimates are quite
sharp: they mainly rely on the covariance structure of $\eta$ and on
Girsanov's theorem applied to $b$.

\vspace{0.3cm}

Before going into the details of our calculations, let us introduce a new set
$\mathcal{\hat{B}}_{M,t}$: as mentioned above, our computations will bring out
some expressions of the form $u_{t}:=\sum_{j\in{\mathbb{Z}}}\delta_{j}\eta
_{j}$, and it will be convenient to keep this kind of term of order
$O(t^{2\alpha-1})$, which is also the order of the exponential correction term
appearing in (\ref{corgirs}). However, since $\delta$ satisfies Proposition
\ref{P:nor_del}, it is easily checked that $u_{t}$ is of the desired order if
$\eta_{j}\leq|j-k|^{\tau}t^{3\alpha-1}$ on $L_{k}$. These considerations
motivate the introduction of the event
\[
\mathcal{\hat{B}}_{M,t}\equiv\{\text{ There exists }\ell\in\bar{\mathbb{Z}%
}_{M}\text{ and }j\in\mathbb{Z}\backslash\{\ell\};\,|\eta_{j}|\geq
|j-\ell|^{\tau}t^{3\alpha-1}\},
\]
and we will trivially bound $\mbox{{\bf P}}(\mathcal{B}_{M,t})$ by
\begin{equation}
\mbox{{\bf P}}(\mathcal{B}_{M,t})\leq\mbox{{\bf P}}(\mathcal{\hat{B}}%
_{M,t})+\mbox{{\bf P}}(\mathcal{\hat{B}}_{M,t}^{c}\cap\mathcal{B}_{M,t}).
\label{tvpatm}%
\end{equation}
We will now prove that the two terms on the right hand side of (\ref{tvpatm})
vanish as $t\rightarrow\infty$, whenever $M$ is large enough.

\vspace{0.5cm}

\noindent\textit{Step 1: Estimation of} $\mbox{{\bf P}}(\mathcal{\hat{B}%
}_{M,t})$

\noindent Let $\Phi$ be the distribution function of a standard Gaussian
random variable, i.e. if $Z\sim\mathcal{N}(0,1)$, then
\begin{equation}
\Phi(x)=\mbox{{\bf P}}(Z\leq x),\label{distgauss}%
\end{equation}
and set $\bar{\Phi}=1-\Phi$. Then let us bound simply
$\mbox{{\bf P}}(\mathcal{\hat{B}}_{M,t})$ by
\begin{align*}
\mbox{{\bf P}}(\mathcal{\hat{B}}_{M,t}) &  \leq\sum_{\ell\in\bar{\mathbb{Z}%
}_{M}}\sum_{j\neq\ell}\mbox{{\bf P}}(|\eta_{j}|\geq|j-\ell|^{\tau}%
t^{3\alpha-1})\\
&  \leq2\sum_{\ell\in\bar{\mathbb{Z}}_{M}}\sum_{j\neq\ell}\bar{\Phi}\left(
\frac{2|j-\ell|^{\tau}t^{\frac{7\alpha-3}{2}}}{C_{0,0}^{1/2}(t)}\right)  ,
\end{align*}
where $C_{0,0}(t)$ defined in (\ref{covet}), equals $t^{\alpha-1}%
/4\mbox{{\bf E}}\left[  \eta_{\ell}\eta_{k}\right]  $. Recall that $\bar{\Phi
}(x)\leq e^{-x^{2}/2}$ for $x$ large enough, and that $C(t)$ satisfies
Proposition \ref{P:ord_t}. Thus, for two constants $c_{1},c_{2}>0$, we get
\begin{equation}
\mbox{{\bf P}}(\mathcal{\hat{B}}_{M,t})\leq c_{1}M\sum_{j\geq1}\exp\left(
-c_{2}j^{2\tau}t^{7\alpha-3}\right)  .\label{bnputm}%
\end{equation}
The following facts are now easily seen:

\begin{itemize}
\item The series in the right hand side of (\ref{bnputm}) is convergent, since
$\tau>0$, which explains the choice of the norm $\Vert x\Vert_{\tau,\ell}$ in
order to bound $\eta_{j}$.

\item Since we have assumed $\alpha>1/2>3/7$, we have $7\alpha-3>0$, and thus,
an elementary application of the dominated convergence theorem yields
\begin{equation}
\lim_{t\rightarrow\infty}\mbox{{\bf P}}(\mathcal{\hat{B}}_{M,t})=0,
\label{limputm}%
\end{equation}
which proves our first claim.
\end{itemize}

\vspace{0.5cm}

\noindent\textit{Step 2: Estimation of} $\mbox{{\bf P}}(\mathcal{\hat{B}%
}_{M,t}^{c}\cap\mathcal{B}_{M,t})$

\noindent Recall that the vector $\delta$ has been introduced because
$-H_{t}(b)-\sum_{j\in{\mathbb{Z}}}\delta_{j}$ $\eta_{j}$ is independent of the
family $\eta$, and for sake of compactness of notations, set
\begin{equation}
e_{t}(b)=\exp\left(  -\beta\left(  H_{t}(b)+\sum_{j\in\mathbb{Z}}\delta
_{j}\eta_{j}\right)  \right)  .\label{defetaeb}%
\end{equation}
Now we have
\begin{align*}
\mbox{{\bf P}}(\mathcal{\hat{B}}_{M,t}^{c}\cap\mathcal{B}_{M,t}) &
=\mbox{{\bf P}}\Bigl(\mathcal{\hat{B}}_{M,t}^{c}\text{ and }E_{b}\left[
\mathbf{1}_{L_{k}}e^{-H_{t}(b)}\right]  <E_{b}\left[  \mathbf{1}_{L_{0}%
}e^{-H_{t}(b)}\right]  \mbox{ for all }k\in\mathbb{\dot{Z}}_{M}\Bigr)\\
&  =\mbox{{\bf P}}\Bigg(\mathcal{\hat{B}}_{M,t}^{c}\text{ and }E_{b}\left[
\mathbf{1}_{L_{k}}e_{t}(b)\exp\Bigl(\sum_{j\in\mathbb{Z}}\beta\delta_{j}%
\eta_{j}\Bigr)\right]  \\
&  \hspace{4cm}<E_{b}\left[  \mathbf{1}_{L_{0}}e_{t}(b)\exp\Bigl(\sum
_{j\in\mathbb{Z}}\beta\delta_{j}\eta_{j}\Bigr)\right]  \mbox{ for all }k\in
\mathbb{\dot{Z}}_{M}\Bigg).
\end{align*}
As mentioned before, $\delta:=C^{-1}(t)v$ depends on the path $b$, as is
easily seen from definition (\ref{defvl}). In order to get rid of the term
$\sum_{j\in\mathbb{Z}}\delta_{j}\eta_{j}$, we will then set
\begin{equation}
\check{\eta}_{0}=\max{(\beta\underline{d}\eta_{0},\beta\overline{d}\eta_{0}%
)},\quad\mbox{ and }\quad\hat{\eta}_{k}=\min{(\beta\underline{d}\eta_{k}%
,\beta\overline{d}\eta_{k})},\label{defchheta}%
\end{equation}
where the constants $\underline{d},\overline{d}$ have been introduced in
Proposition \ref{P:nor_del}. Then, according to the definition of
$\mathcal{\hat{B}}_{M,t}^{c}$, we get
\begin{align*}
&  \mbox{{\bf P}}(\mathcal{\hat{B}}_{M,t}^{c}\cap\mathcal{B}_{M,t})\\
\leq &  \mbox{{\bf P}}\Bigg(\mbox{For any }k\in\mathbb{\dot{Z}}_{M}%
,\,E_{b}\left[  \mathbf{1}_{L_{k}}e_{t}(b)\exp\Bigl(-\sum_{j\in\mathbb{Z}%
}\beta|\delta_{j}||j-k|^{\tau}t^{3\alpha-1}+\hat{\eta}_{k}\Bigr)\right]  \\
&  \hspace{5cm}<E_{b}\left[  \mathbf{1}_{L_{0}}e_{t}(b)\exp\Bigl(\sum
_{j\in\mathbb{Z}}\beta|\delta_{j}|j^{\tau}t^{3\alpha-1}+\check{\eta}%
_{0}\Bigr)\right]  \Bigg).
\end{align*}
Now, invoking Proposition~\ref{P:nor_del} item (iii), we obtain that for any
integer $k$, there exists a constant $\gamma$ (possibly depending on $\beta$)
such that $\sum_{j\in\mathbb{Z}}\beta|\delta_{j}||j-\ell|^{\tau}\leq\gamma
t^{-\alpha}$ on $L_{k}$. Thus, thanks to the fact that the random variables
$\eta$ only depend on $W$, and observing that $\tilde{Z}_{t}^{\alpha}%
(k)=E_{b}[\mathbf{1}_{L_{k}}e_{t}(b)]$, we get
\begin{align}
&  \mbox{{\bf P}}(\mathcal{\hat{B}}_{M,t}^{c}\cap\mathcal{B}_{M,t})\nonumber\\
\leq &  \mbox{{\bf P}}\Big(\mbox{For any }k\in\mathbb{\dot{Z}}_{M};\,\tilde
{Z}_{t}^{\alpha}(k)\exp(-\gamma t^{2\alpha-1}+\hat{\eta}_{k})<\exp(\gamma
t^{2\alpha-1}+\check{\eta}_{0})\tilde{Z}_{t}^{\alpha}(0)\Bigr)\nonumber\\
= &  \mbox{{\bf P}}\Biggl(\mbox{For any }k\in\mathbb{\dot{Z}}_{M}%
;\,\frac{\tilde{Z}_{t}^{\alpha}(k)}{\tilde{Z}_{t}^{\alpha}(0)}<\exp(2\gamma
t^{2\alpha-1}+\check{\eta}_{0}-\hat{\eta}_{k})\Biggr).\label{whenQcondchanges}%
\end{align}
Let us apply now Proposition~\ref{P:entr} in order to conclude that
\[
\mbox{{\bf P}}(\mathcal{\hat{B}}_{M,t}^{c}\cap\mathcal{B}_{M,t})\leq
\mbox{{\bf P}}\Biggl(\mbox{For any }k\in\mathbb{\dot{Z}}_{M};\,\frac{\tilde
{Z}_{t}^{\alpha}(0,W^{k,t})}{\tilde{Z}_{t}^{\alpha}(0,W)}<\exp(\hat{\gamma
}t^{2\alpha-1}+\check{\eta}_{0}-\hat{\eta}_{k})\Biggr),
\]
where $\hat{\gamma}=\hat{\gamma}(M)=\sup\{2\gamma+\zeta(k);k\in\mathbb{\dot
{Z}}_{M}\}$ and $\zeta\left(  k\right)  =4k(k+1)$. We have thus proved that
\[
\mbox{{\bf P}}(\mathcal{\hat{B}}_{M,t}^{c}\cap\mathcal{B}_{M,t})\leq
\mbox{{\bf P}}(F_{M}),
\]
where
\[
F_{M}=\left\{  \mbox{For any }k\in\mathbb{\dot{Z}}_{M};\,\frac{\tilde{Z}%
_{t}^{\alpha}(0,W^{k,t})}{\tilde{Z}_{t}^{\alpha}(0,W)}<\exp(\hat{\gamma
}t^{2\alpha-1}+\check{\eta}_{0}-\hat{\eta}_{k})\right\}  .
\]

\vspace{0.5cm}

\noindent\textit{Step 3: Evaluation of} $\mbox{{\bf P}}(F_{M})$

\noindent We can see now that the probability of $F_{M}$ will be expressed in
terms of a balance between the values of $\check{\eta}_{0}-\hat{\eta}_{k}$
(which will be assumed to be highly negative) and the ratio $\tilde{Z}%
_{t}^{\alpha}(0,W^{k,t})/\tilde{Z}_{t}^{\alpha}(0,W)$, which cannot be too
small at many different sites $k$. In order to quantify this heuristic
statement, we introduce a family $\mathcal{\bar{S}}_{M,m}$ of subsets of
${\bar{\mathbb{Z}}}_{M}$ which will be used to construct a large symmetric set
$L$ around 0 such that $\check{\eta}_{0}-\hat{\eta}_{\ell}<-t^{2\alpha-1+\rho
}$ for all $\ell\in L$:
for a given $\rho>0$ and integer numbers $m$ and $M$, define the families of
subsets
\begin{align}
\mathcal{S}_{M,m}  &  =\bigcup_{k,\hat{k}\in D_{M,m}}\left\{  k{\mathbb{\dot
{Z}}}_{\hat{k}}\right\}  ,\quad\mbox{ with }\quad D_{M,m}=\left\{  \left(
k,k^{\prime}\right)  :k\geq1,\hat{k}\geq m;k{\mathbb{\dot{Z}}}_{\hat{k}%
}\subset\mathbb{\dot{Z}}_{M}\right\} \nonumber\\
\mathcal{\bar{S}}_{M,m}  &  =\left\{  L\subset\mathbb{\dot{Z}}_{M}%
;\mbox{ There exists }S\in\mathcal{S}_{M,m}\mbox{ such that }S\subset
L\right\}  . \label{defbsmm}%
\end{align}
In relation with these families of subsets of $\mathbb{\dot{Z}}_{M}$, set
also
\begin{equation}
\hat{F}_{M,m,\rho}=\bigcup_{L\in\mathcal{\bar{S}}_{M,m}}\hat{F}_{\rho,L},
\label{defhfmr}%
\end{equation}
with
\begin{equation}
\hat{F}_{\rho,L}=\big\{\check{\eta}_{0}-\hat{\eta}_{\ell}<-t^{2\alpha-1+\rho
},\mbox{ for all }\ell\in L,\check{\eta}_{0}-\hat{\eta}_{\hat{\ell}%
}>-t^{2\alpha-1+\rho},\mbox{ for all }\hat{\ell}\in\mathbb{\dot{Z}}%
_{M}\backslash L\big\}. \label{defhfmr2}%
\end{equation}
Then one can bound trivially $\mbox{{\bf P}}(F_{M})$ by
\[
\mbox{{\bf P}}(F_{M})\leq1-\mbox{{\bf P}}(\hat{F}_{M,m,\rho}%
)+\mbox{{\bf P}}(F_{M}\cap\hat{F}_{M,m,\rho}).
\]
Furthermore, for $t$ large enough, we have $\hat{\gamma}t^{2\alpha
-1}-t^{2\alpha-1+\rho}<0$, which explains the need for the constant $\rho>0$.
Thus%
\begin{align*}
F_{M}\cap\hat{F}_{M,m,\rho}  &  \subseteq\bigcup_{L\in\mathcal{\bar{S}}_{M,m}%
}\bigcap_{\ell\in L}\left\{  \frac{\tilde{Z}_{t}^{\alpha}(0,W^{\ell,t}%
)}{\tilde{Z}_{t}^{\alpha}(0,W)}<\exp\left(  \hat{\gamma}t^{2\alpha
-1}-2t^{2\alpha-1+\rho}\right)  \right\}  \cap\hat{F}_{\rho,L}\\
&  \subseteq\bigcup_{L\in\mathcal{\bar{S}}_{M,m}}\left\{  \tilde{Z}%
_{t}^{\alpha}(0,W^{\ell,t})<\tilde{Z}_{t}^{\alpha}(0,W)\mbox{ for
all }\ell\in L\right\}  \cap\hat{F}_{\rho,L}.
\end{align*}
Hence, we get
\begin{align}
\mbox{{\bf P}}(F_{M})  &  \leq1-\mbox{{\bf P}}(\hat{F}_{M,m,\rho})+\sum
_{L\in\mathcal{\bar{S}}_{M,m}}\mbox{{\bf P}}\left(  \left\{  \tilde{Z}%
_{t}^{\alpha}(0,W^{\ell,t})<\tilde{Z}_{t}^{\alpha}(0,W)\mbox{ for
all }\ell\in L\right\}  \cap\hat{F}_{\rho,L}\right) \nonumber\\
&  \leq1-\mbox{{\bf P}}(\hat{F}_{M,m,\rho})+\sum_{L\in\mathcal{\bar{S}}_{M,m}%
}\mbox{{\bf P}}\left(  \tilde{Z}_{t}^{\alpha}(0,W^{\ell,t})<\tilde{Z}%
_{t}^{\alpha}(0,W)\mbox{ for
all }\ell\in L\right)  \mbox{{\bf P}}\left(  \hat{F}_{\rho,L}\right)  ,
\label{E:fin}%
\end{align}
where in the last step, we have used the independence, proved in the next
step, between the random variables $\tilde{Z}_{t}^{\alpha}(0,W^{\ell,t})$ and
the sequence $\{\eta_{k};k\in{\bar{\mathbb{Z}}}_{M}\}$, and also between
$\tilde{Z}_{t}^{\alpha}(0,W)$ and the sequence $\{\eta_{k};k\in{\bar
{\mathbb{Z}}}_{M}\}$.\vspace{0.5cm}

\noindent\emph{Step 4: Independence of }$\eta$\emph{ and the }$\tilde{Z}%
_{t}^{\alpha}$\emph{'s.}

\noindent Using the notation $X\left(  W,b\right)  $ introduced in Lemma
\ref{independence}, this lemma's conclusion is that $X\left(  W,b\right)  $
and $\eta\left(  W\right)  $ are independent for each continuous function $b$;
after evaluation of $\tilde{Z}_{t}^{\alpha}(0,W)$ in formula (\ref{Ztildadef}%
), it implies that the latter is also independent of $\eta$.

Lemma \ref{independence} can also be applied to prove the other independence:
it proves that for each fixed $b,k$, we have independence of $X\left(
W^{k,t},b\right)  $ and the entire sequence $\eta$. When defining $\tilde
{Z}_{t}^{\alpha}(0,W^{\ell,t})$, formula (\ref{Ztildadef}) must be used with
$W$ replaced by $W^{\ell,t}$, which specifically means%
\begin{align*}
\tilde{Z}_{t}^{\alpha}(0,W^{\ell,t}) &  =E_{b}\left[  \mathbf{1}_{L_{k}}%
\exp\beta\left(  \int_{0}^{t}W^{\ell,t}(ds,b_{s})-\sum_{j\in\mathbb{Z}}%
\delta_{j}\left(  b,\mathcal{L}\left(  W^{\ell,t}\right)  \right)  \eta
_{j}\left(  W^{\ell,t}\right)  \right)  \right]  \\
&  =E_{b}\left[  \mathbf{1}_{L_{k}}\exp\beta\left(  X\left(  W^{\ell
,t},b\right)  \right)  \right]  ,
\end{align*}
proving that $\tilde{Z}_{t}^{\alpha}(0,W^{\ell,t})$ is independent of $\eta$,
as required to justify (\ref{E:fin}) in step 3.

One can prove in addition that $\delta\left(  b,\mathcal{L}\left(  W\right)
\right)  =\delta\left(  b,\mathcal{L}\left(  W^{\ell,t}\right)  \right)  $ for
any $\ell$, but this fact will not be needed.

\vspace{0.3cm}

\noindent\emph{Step 5: finishing the proof. }The end of our proof of Lemma
\ref{L:env} relies on the following propositions, whose proofs will be
postponed until the next sections.

\begin{proposition}
\label{P:des} Let $m$ be a fixed positive even integer, and $M>m$. Then, for
any $L\in\mathcal{\bar{S}}_{M,m}$, we have
\[
\mbox{{\bf P}}\left(  \tilde{Z}_{t}^{\alpha}(0,W^{\ell,t})<\tilde{Z}%
_{t}^{\alpha}(0,W)\mbox{ for all }\ell\in L\right)  \leq\frac{1}{m}.
\]

\end{proposition}

\begin{proposition}
\label{P:troi} Let $m$ be a fixed positive integer. Let $\rho$ be a strictly
positive number such that $\frac{5}{2}(\alpha-\frac{3}{5})$ $+\rho<0$. Then,
for $t$ large enough, there exists a $M$ large enough such that
\begin{equation}
\mbox{{\bf P}}(\hat{F}_{M,m,\rho})\geq1-\frac{1}{m}. \label{majpfmr}%
\end{equation}

\end{proposition}

With these results in mind, let us finish now the proof of Lemma \ref{L:env},
and thus of our theorem: take $t,M$ large enough so that (\ref{majpfmr}) is
satisfied. Then (\ref{E:fin}) yields directly, invoking Proposition
\ref{P:des} and the fact that the events $\hat{F}_{\rho,L}$ are disjoints,
\[
\mbox{{\bf P}}(F_{M})\leq\frac{1}{m}+\frac{1}{m}\sum_{L\in\mathcal{\bar{S}%
}_{M,m}}\mbox{{\bf P}}(\hat{F}_{\rho,L})\leq\frac{1}{m}+\frac{1}{m}=\frac
{2}{m}.
\]
which tends to 0 as $m\rightarrow\infty$, and ends the proof of the theorem,
modulo establishing the last two propositions above.\hfill$\Box$

\vspace{0.3cm}

Before proceeding with the proofs of Propositions \ref{P:des} and
\ref{P:troi}, we discuss the consequences of weakening Hypothesis \ref{env1}.
If we assume only that%
\begin{equation}
Q\left(  x\right)  \leq\left\vert x\right\vert ^{-2-\theta}, \label{newQcond}%
\end{equation}
can we find values of $\theta\leq1$ such that we still get superdiffusive
behavior for the polymer, i.e. $\alpha>1/2$? Since the result of the Girsanov
theorem, Proposition \ref{P:entr}, is not effected by the value of $\theta$
above, this means that the penalization from Girsanov's theorem, of order
$t^{2\alpha-1}$, cannot be made smaller by a different choice of decorrelation
speed in $Q$. Therefore we should expect not to be able to preserve the
threshold $\alpha<3/5$. To see exactly what happens to this threshold under
condition (\ref{newQcond}), we first state, and leave it to the reader to
check, that we can rework the proof of Proposition~\ref{P:nor_del} item (iii)
to obtain instead
\[
\left\vert \delta\right\vert _{\tau,k}-\delta_{k}=o\left(  t^{-\alpha\theta
}\right)  .
\]

It is then simple to check that (\ref{whenQcondchanges}) becomes%
\[
\mbox{{\bf P}}(\mathcal{\hat{B}}_{M,t}^{c}\cap\mathcal{B}_{M,t})\leq
\mbox{{\bf P}}\Biggl(\mbox{For any }k\in\mathbb{\dot{Z}}_{M};\,\frac{\tilde
{Z}_{t}^{\alpha}(k)}{\tilde{Z}_{t}^{\alpha}(0)}<\exp(2\gamma t^{3\alpha
-1-\theta}+\check{\eta}_{0}-\hat{\eta}_{k})\Biggr).
\]
Hence the application of Proposition \ref{P:entr} still works, but we can no
longer make the corresponding Girsanov penalization of the same order, since
for $\theta<1$, $3\alpha-1-\alpha\theta>2\alpha-1$. Having thus convinced
ourselves that Hypothesis \ref{env1} is the only way to get the entire proof
to be efficient in terms of using comparable penalizations throughout, we can
now ignore this inefficiency, and answer the question at the beginning of this
paragraph. The reader will check that any other occurences of the use of
Hypothesis \ref{env1} are not further effected by switching to (\ref{newQcond}%
): the entire proof can still be used if we only require that the magnitude of
the $\eta_{k}$'s, namely $t^{\left(  1-\alpha\right)  /2}$, is larger than the
new penalization $t^{3\alpha-1-\alpha\theta}$. This yields%
\[
\alpha<\frac{3}{7-2\theta}.
\]
Now we see that to get a super-diffusive behavior, we need $3/\left(
7-2\theta\right)  >1/2$, i.e. $\theta>1/2$. We also see that the weakest
hypothesis required for such behavior is $Q\left(  x\right)  \leq
x^{-5/2-\vartheta}$ for $\vartheta>0.$ We state these findings formally, using
the reparametrization $\theta=\vartheta+1/2$.

\begin{corollary}
\label{disorder}Assume instead of Hypothesis \ref{env1} that there exists
$\vartheta\in(0,1/2]$ such that as $\left\vert x\right\vert \rightarrow\infty
$,%
\[
Q\left(  x\right)  =O\left(  \left\vert x\right\vert ^{-5/2-\vartheta}\right)
.
\]
Then for any $\varepsilon>0$ we obtain the following specific super-diffusive
behavior for the polymer measure:%
\[
\lim_{t\rightarrow\infty}\mathbf{P}\left[  \langle\sup_{s\leq t}|b_{s}%
|\rangle_{t}\geq t^{\frac{1}{2}+\frac{\vartheta}{6-2\vartheta}-\varepsilon
}\right]  =1.
\]

\end{corollary}

\subsection{Proof of Proposition \ref{P:des}}

Let $L\in\mathcal{\bar{S}}_{M,m}$. Then, by definition (\ref{defbsmm}) of
$\mathcal{\bar{S}}_{M,m}$, there exists $k\geq1$ such that $k{\mathbb{\dot{Z}%
}}_{m}\subset L$. Then
\begin{align*}
&  \mbox{{\bf P}}\left(  \tilde{Z}_{t}^{\alpha}(0,W^{\ell,t})<\tilde{Z}%
_{t}^{\alpha}(0,W)\mbox{ for all }\ell\in L\right) \\
\leq &  \mbox{{\bf P}}\left(  \tilde{Z}_{t}^{\alpha}(0,W^{\ell,t})<\tilde
{Z}_{t}^{\alpha}(0,W)\mbox{ for all }\ell\in k{\mathbb{\dot{Z}}}_{m}\right)  .
\end{align*}
It is thus sufficient to estimate the right hand side in the above
inequality.\vspace{0.3cm}

Given an even integer $m\leq M$, recall that ${\bar{\mathbb{Z}}}_{m}$ has been
defined at (\ref{defbzm}). Set also $\hat{m}=m/2$, and for each $i\in
k{\bar{\mathbb{Z}}}_{\hat{m}}$, we associate the following event:
\[
\Omega^{(i)}\equiv\left\{  \tilde{Z}_{t}^{\alpha}(0,W^{\ell,t})<\tilde{Z}%
_{t}^{\alpha}(0,W^{i,t})\mbox{ for all }\ell\in k{\bar{\mathbb{Z}}}_{\hat{m}%
}\backslash\{i\}\right\}  .
\]
Then these events are disjoint, and since $|k{\bar{\mathbb{Z}}}_{\hat{m}%
}|=2\hat{m}+1$, we get trivially the existence of $i_{0}\in k{\bar{\mathbb{Z}%
}}_{\hat{m}}$ such that
\begin{equation}
\mbox{{\bf P}}\left(  \Omega^{(i_{0})}\right)  \leq\frac{1}{2\hat{m}+1}%
\leq\frac{1}{m}.\label{inqoio}%
\end{equation}
However, the translation-invariance of the environment $W$ yields
\begin{align}
&  \mbox{{\bf P}}\left(  \tilde{Z}_{t}^{\alpha}(0,W^{\ell,t})<\tilde{Z}%
_{t}^{\alpha}(0,W)\mbox{ for all }\ell\in k{\mathbb{\dot{Z}}}_{m}\right)
\nonumber\\
&  =\mbox{{\bf P}}\left(  \tilde{Z}_{t}^{\alpha}(0,W^{\ell+i_{0},t})<\tilde
{Z}_{t}^{\alpha}(0,W^{i_{0},t})\mbox{ for all }\ell\in k{\mathbb{\dot{Z}}}%
_{m}\right)  .\nonumber
\end{align}
Indeed, exactly as we proved Lemma \ref{WWkt}, denoting again $h\left(
s\right)  =\min(2s/t,1)2kt^{\alpha}$, it holds that for fixed $b$, $\int
_{0}^{t}W\left(  ds,b_{s}+\left(  \ell+i_{0}\right)  h\left(  s\right)
\right)  $ has the same distribution as $\int_{0}^{t}W(ds,b_{s}+$ $\ell
h\left(  s\right)  )$.

We may now rewrite the above expression as the following upper bound:%
\begin{align}
&  \mbox{{\bf P}}\left(  \tilde{Z}_{t}^{\alpha}(0,W^{\ell,t})<\tilde{Z}%
_{t}^{\alpha}(0,W)\mbox{ for all }\ell\in k{\mathbb{\dot{Z}}}_{m}\right)
\nonumber\\
&  \leq\mbox{{\bf P}}\left(  \tilde{Z}_{t}^{\alpha}(0,W^{\ell,t})<\tilde
{Z}_{t}^{\alpha}(0,W^{i_{0},t})\mbox{ for all }\ell\in k{\bar{\mathbb{Z}}%
}_{\hat{m}}\backslash\{i_{0}\}\right)  =\mbox{{\bf P}}\left(  \Omega^{(i_{0}%
)}\right)  .\label{inqtz}%
\end{align}
Observe that the last inequality is just due to the elementary fact that
$k{\bar{\mathbb{Z}}}_{\hat{m}}\backslash\{i_{0}\}\subset i_{0}+k$%
{$\mathbb{\dot{Z}}$}$_{m}$ whenever $i_{0}\in k{\bar{\mathbb{Z}}}_{\hat{m}}$,
a fact which is easily checked.
Hence, putting together (\ref{inqoio}) and (\ref{inqtz}), we get the announced
result.\hfill$\Box$

\subsection{Proof of Proposition \ref{P:troi}}

Recall that $\hat{F}_{M,m,\rho}$ is defined by (\ref{defhfmr}), and define the
quantity
\[
\tau(t):=2\beta^{-1}t^{\frac{5}{2}(\alpha-\frac{3}{5})+\rho},
\]
which tends to 0 as $t\rightarrow\infty$ if $\alpha<\frac{3}{5}$ and $\rho$ is
small enough. The following inequality
\begin{equation}
\mbox{{\bf P}}(\hat{F}_{M,m,\rho})\geq\mbox{{\bf P}}\left(  \bigcup
_{L\in\mathcal{\bar{S}}_{M,m}}\left\{  t^{(\alpha-1)/2}(\check{\eta}_{0}%
-\hat{\eta}_{\ell})\leq-\beta\tau(t)\mbox{ for all }\ell\in L\right\}
\right)  \label{lbpfhat}%
\end{equation}
is then easily established by an elementary inclusion argument, which we
detail here. Indeed, assume that for some $L\in\mathcal{\bar{S}}_{M,m}$, for
all $\ell\in L$, $\eta$ satisfies%
\[
t^{(\alpha-1)/2}(\check{\eta}_{0}-\hat{\eta}_{\ell})\leq-\beta\tau(t)
\]
which is equivalent to%
\[
\check{\eta}_{0}-\hat{\eta}_{\ell}\leq-t^{2\alpha-1+\rho}.
\]
To justify the above inequality, we only need to prove that for some other
$L^{\prime}\in\mathcal{\bar{S}}_{M,m}$, the same $\eta$ also satisfies the
above inequality for all $\ell\in L^{\prime}$, while for all $\ell
\in\mathbb{\dot{Z}}_{M}\setminus L^{\prime}$, the contrary holds, namely%
\[
\check{\eta}_{0}-\hat{\eta}_{\ell}>-t^{2\alpha-1+\rho}.
\]
Let then $\Lambda$ be the subset of $\mathbb{\dot{Z}}_{M}$ defined by
\[
\Lambda=\left\{  \ell\in\mathbb{\dot{Z}}_{M};\,\check{\eta}_{0}-\hat{\eta
}_{\ell}>-t^{2\alpha-1+\rho}\right\}  ,
\]
and set $L^{\prime}=\mathbb{\dot{Z}}_{M}\setminus\Lambda$. Then, by
construction $L^{\prime}$ has the required properties defined above, and since
$L^{\prime}\supset L$, by definition of $\mathcal{\bar{S}}_{M,m}$, we have
$L^{\prime}\in\mathcal{\bar{S}}_{M,m}$.

In order to get a lower bound on the right hand side of (\ref{lbpfhat}) above,
we will construct now a large enough collection of symmetric and disjoint sets
in $\mathbb{\dot{Z}}_{M}$: with $m<M$, consider the collection $\{Q_{q}%
(m)${$\mathbb{\dot{Z}}$}$_{m};q<q^{\ast}\}$, where the integers $Q_{q}\left(
m\right)  $ are defined by
\[
Q_{1}(m)=1,\hspace{0.5cm}Q_{q+1}(m)=mQ_{q}(m)+1,\hspace{0.5cm}q^{\ast}%
=\inf{\{q;Q_{q}(m)>M\}}.
\]
This collection is the sequence%
\[
{\mathbb{\dot{Z}}}_{m},(m+1){\mathbb{\dot{Z}}}_{m},[m(m+1)+1]{\mathbb{\dot{Z}%
}}_{m},\cdots,Q_{q}\left(  m\right)  {\mathbb{\dot{Z}}}_{m},\cdots,Q_{q^{\ast
}-1}\left(  m\right)  {\mathbb{\dot{Z}}}_{m},
\]
which are non-overlapping annuli in $\mathbb{\dot{Z}}_{M}$, and therefore are
indeed symmetric and disjoint subsets of $\mathbb{\dot{Z}}_{M}$.
Since $Q_{q}\left(  m\right)  ${$\mathbb{\dot{Z}}$}$_{m}$ is certainly of the
form $k{\mathbb{\dot{Z}}}_{\hat{k}}$ with $k\geq1$ and $\hat{k}\geq m$, and is
a subset of $\mathbb{\dot{Z}}_{M}$ as soon as $q<q^{\ast}$, by definition
$Q_{q}\left(  m\right)  ${$\mathbb{\dot{Z}}$}$_{m}\in\mathcal{\bar{S}}_{M,m}$.
Thus, using the notation $\check{\eta}_{0},\hat{\eta}_{\ell}$ and $\eta_{\ell
}$ defined in (\ref{defetaeb}) and (\ref{defchheta}), and reverting to the
notation $\tilde{\eta}=2t^{-\left(  1-\alpha\right)  /2}\eta$, we get%
\[
\mbox{{\bf P}}(\hat{F}_{M,m,\rho})\geq\mbox{{\bf P}}\Bigg(\bigcup_{q<q^{\ast}%
}\big\{\max(\underline{d}\tilde{\eta}_{0},\overline{d}\tilde{\eta}_{0}%
)-\min(\underline{d}\tilde{\eta}_{\ell},\overline{d}\tilde{\eta}_{\ell}%
)\leq-\tau(t)\mbox{ for all }\ell\in Q_{q}(m){\mathbb{\dot{Z}}}_{m}%
\big\}\Bigg).
\]
Indeed, the original set $\hat{F}_{M,m,\rho}$ defined in (\ref{defhfmr}),
(\ref{defhfmr2}) was a union of events indexed by $L\in\mathcal{\bar{S}}%
_{M,m}$, while here we use only sets of the form $L=Q_{q}\left(  m\right)
${$\mathbb{\dot{Z}}$}$_{m}$; moreover, the above condition on the difference
$\max(\underline{d}\tilde{\eta}_{0},\overline{d}\tilde{\eta}_{0}%
)-\min(\underline{d}\tilde{\eta}_{\ell},\overline{d}\tilde{\eta}_{\ell})$ is
implied by the two conditions on the individual terms of this difference in
$\hat{F}_{M,m,\rho}$, and the shorthand notation $\tau\left(  t\right)  $ was
introduced above to be consistent with these conditions in (\ref{defhfmr2}).
Let us call now $A_{\ell}$ the event
\[
A_{\ell}=\left\{  \max(\underline{d}\tilde{\eta}_{0},\overline{d}\tilde{\eta
}_{0})-\min(\underline{d}\tilde{\eta}_{\ell},\overline{d}\tilde{\eta}_{\ell
})\leq-\tau(t)\right\}  ,
\]
and we distinguish two cases according to the values of $\tilde{\eta}_{0}$:

\begin{itemize}
\item[\textbf{(a)}] If $\tilde{\eta}_{0}\geq0$, then $\max(\underline{d}%
\tilde{\eta}_{0},\overline{d}\tilde{\eta}_{0})=\overline{d}\tilde{\eta}_{0}$,
and hence $A_{\ell}$ is the event defined by the relation
\[
\min(\underline{d}\tilde{\eta}_{\ell},\overline{d}\tilde{\eta}_{\ell})\geq
\tau(t)+\overline{d}\tilde{\eta}_{0}.
\]
In particular, $\tilde{\eta}_{\ell}$ has to be positive, and thus $A_{\ell}$
can be written as
\[
\left\{  \overline{d}\tilde{\eta}_{0}-\underline{d} \tilde{\eta}_{\ell}%
<-\tau(t)\right\}  .
\]

\item[\textbf{(b)}] If $\tilde{\eta}_{0}\leq-\tau(t)/\underline{d}\leq0$, then
$\max(\underline{d}\tilde{\eta}_{0},\overline{d}\tilde{\eta}_{0}%
)=\underline{d}\tilde{\eta}_{0}$. Thus $A_{\ell}$ can be written as the event
defined by the relation%
\begin{equation}
\min(\underline{d}\tilde{\eta}_{\ell},\overline{d}\tilde{\eta}_{\ell})\geq
\tau(t)+\underline{d}\tilde{\eta}_{0},\label{aleq}%
\end{equation}
and if $\tilde{\eta}_{0}\leq-\tau(t)/\underline{d}$, we have $\tau
(t)+\underline{d}\tilde{\eta}_{0}\leq0$. Hence, (\ref{aleq}) is implied by
$\tilde{\eta}_{\ell}\geq0$.
\end{itemize}

Summarizing the considerations above, we get
\[
\mbox{{\bf P}}(\hat{F}_{M,m,\rho})\geq\mbox{{\bf P}}(D^{+}%
)+\mbox{{\bf P}}(D^{-}),
\]
with
\begin{align*}
D^{+}  &  =\bigcup_{q<q^{\ast}}\left\{  \overline{d}\tilde{\eta}%
_{0}-\underline{d}\tilde{\eta}_{\ell}\leq-\tau(t)\mbox{ for all
}\ell\in Q_{q}(m){\mathbb{\dot{Z}}}_{m}\right\}  \cap\left\{  \tilde{\eta}%
_{0}>0\right\} \\
D^{-}  &  =\bigcup_{q<q^{\ast}}\left\{  \tilde{\eta}_{\ell}\geq
0\mbox{ for all
}\ell\in Q_{q}(m){\mathbb{\dot{Z}}}_{m}\right\}  \cap\left\{  \tilde{\eta}%
_{0}\leq-\tau(t)/\underline{d}\right\}  .
\end{align*}
We will now prove that $\mbox{{\bf P}}(D^{+})$ is close to $1/2$. Entirely
similar arguments, left to the reader, lead to showing that
$\mbox{{\bf P}}(D^{-})$ can also be made arbitrarily close to $1/2$,
concluding the proof of the proposition.\vspace{0.3cm}

Observe that, according to Proposition \ref{P:ord_t} the random variables
$\{\tilde{\eta}_{\ell};l\in{\bar{{\mathbb{Z}}}}_{M}\}$ converge in
distribution to a family of independent standard Gaussian random variables
$\{\Upsilon_{\ell};l\in{\bar{{\mathbb{Z}}}}_{M}\}$. Consequently, and using
the fact that $-\tau\left(  t\right)  \rightarrow0$ as $t\rightarrow\infty$,
\[
\mbox{{\bf P}}(D^{+})=\mbox{{\bf P}}\left(  \bigcup_{q<q^{\ast}}\left\{
\overline{d}\Upsilon_{0}-\underline{d}\Upsilon_{\ell}\leq0\mbox{ for all }\ell
\in Q_{q}(m){\mathbb{\dot{Z}}}_{m}\right\}  \cap\left\{  \Upsilon
_{0}>0\right\}  \right)  +\varepsilon_{M}(t),
\]
where, for a fixed $M\in{\mathbb{N}}$, we have $\lim_{t\rightarrow\infty
}\varepsilon_{M}(t)=0$. Furthermore, since the $\Upsilon_{\ell}$ are
independent random variables, we get
\begin{align}
\mbox{{\bf P}}(D^{+}) &  =\int_{0}^{\infty}\mbox{{\bf P}}\Bigg(\bigcup
_{q<q^{\ast}}\left\{  \overline{d}x-\underline{d}\Upsilon_{\ell}%
\leq0\mbox{ for all }\ell\in Q_{q}(m){\mathbb{\dot{Z}}}_{m}\right\}
\Bigg)\frac{e^{-\frac{x^{2}}{2}}}{(2\pi)^{1/2}}dx+\varepsilon_{M}%
(t)\nonumber\\
&  =\frac{1}{2}-\int_{0}^{\infty}\mbox{{\bf P}}\Bigg(\bigcap_{q<q^{\ast}}%
\hat{D}_{q}\Bigg)\frac{e^{-\frac{x^{2}}{2}}}{(2\pi)^{1/2}}dx+\varepsilon
_{M}(t),\label{inqdp}%
\end{align}
where
\[
\hat{D}_{q}=\left\{  \mbox{There exists }\ell\in Q_{q}(m){\mathbb{\dot{Z}}%
}_{m};\,\overline{d}x-\underline{d}\Upsilon_{\ell}\geq0\right\}  .
\]
In order to take advantage of the independence of the $\Upsilon_{\ell}$, it is
convenient to pick some disjoint sets out of $\mathbb{\dot{Z}}_{M}$, which
explains the choice of disjoint subsets $Q_{q}(m)${$\mathbb{\dot{Z}}$}$_{m}$.
Now, it is easily seen that, for a fixed value $q_{0}$, if one desires to have
$q^{\ast}>q_{0}$, it is sufficient to take $M$ of order $m^{q_{0}}$. Let us
assume that we are in this situation; this means that, setting $\kappa
=\overline{d}/\underline{d}$, we have
\begin{align*}
\mbox{{\bf P}}\Bigg(\bigcap_{q<q^{\ast}}\hat{D}_{q}\Bigg) &  \leq
\mbox{{\bf P}}\Bigg(\bigcap_{q\leq q_{0}}\left\{  \mbox{There exists }\ell\in
Q_{q}(m){\mathbb{\dot{Z}}}_{m};\,\Upsilon_{\ell}\leq\kappa x\right\}  \Bigg)\\
&  =\mbox{{\bf P}}^{q_{0}}\left(  \mbox{There exists }\ell\in{\mathbb{\dot{Z}%
}}_{m};\,\Upsilon_{\ell}\leq\kappa x\right)  =\left[  1-\mbox{{\bf P}}^{2m}%
\left(  \Upsilon_{1}\geq\kappa x\right)  \right]  ^{q_{0}}.
\end{align*}
Plugging these inequalities into (\ref{inqdp}), we obtain
\[
\mbox{{\bf P}}(D^{+})\geq\frac{1}{2}-\int_{0}^{\infty}\left[
1-\mbox{{\bf P}}^{2m}\left(  \Upsilon_{1}\geq\kappa x\right)  \right]
^{q_{0}}\frac{e^{-\frac{x^{2}}{2}}}{(2\pi)^{1/2}}dx+\varepsilon_{M}\left(
t\right)  .
\]
Recall that the functions $\Phi$ has been defined by relation (\ref{distgauss}%
). Then the last inequality yields,
\[
\mbox{{\bf P}}(D^{+})\geq\frac{1}{2}-\int_{0}^{\infty}\left[  1-\Phi\left(
\kappa x\right)  ^{2m}\right]  ^{q_{0}}\frac{e^{-\frac{x^{2}}{2}}}%
{(2\pi)^{1/2}}dx+\varepsilon_{M}\left(  t\right)  .
\]
It is easily seen that this probability can be made as close as we wish to
$\frac{1}{2}$ by taking $q_{0}\rightarrow\infty$, because $1/2\leq\Phi(x)<1$
for all $x\geq0$, this asymptotic being equivalent to $M\rightarrow\infty
$.\hfill$\Box$

\end{document}